\documentclass[reqno,a4paper,11pt]{amsart}
 
\usepackage{amssymb,a4,amsthm,amsmath,enumerate}

\usepackage{texdraw}
\everytexdraw{\drawdim pt \linewd 0.5 \textref h:C v:C \setunitscale 1}

\usepackage[all]{xy}
\newcommand{\xycenter}[1]{\vcenter{\hbox{\xymatrix{#1}}}}
\SelectTips{cm}{}
\newdir{ >}{{}*!/-7pt/@{>}}

\newcommand{\dlpullback}[1][dl]{\save*!/#1-1.2pc/#1:(-1,1)@^{|-}\restore}

\newcommand{\drpullback}[1][dr]{\save*!/#1-1.2pc/#1:(-1,1)@^{|-}\restore}

\usepackage{mathrsfs}

\swapnumbers

\newtheoremstyle{myprop}%
{3pt}
{}
{\itshape}
{}
{\bf}
{.}
{.5em}
{}%

\newtheoremstyle{mydef}%
{3pt}
{3pt}
{}
{}
{\itshape}
{.}
{.5em}
{}%

\theoremstyle{myprop}
\newtheorem{proposition}{Proposition}[section]
\newtheorem{corollary}[proposition]{Corollary}
\newtheorem{theorem}[proposition]{Theorem}
\newtheorem{lemma}[proposition]{Lemma}

\theoremstyle{mydef}

\newtheorem{examples}[proposition]{Examples}
\newtheorem{example}[proposition]{Example}

\theoremstyle{remark}

\newtheorem{para}[proposition]{} 

\newcommand{\myemph}{\textit}

\newcommand{\Nat}{\mathbb{N}}
\DeclareMathOperator*{\colim}{colim}
\newcommand{\comp}{\,}
\newcommand{\bcomp}{\circ}
\newcommand{\iso}{\cong}

\newcommand{\adjoint}{\dashv}
\newcommand{\mysup}{\mathrm{sup}}

\newcommand{\defcat}{\mathrm}

\newcommand{\Set}{\defcat{Set}}
\newcommand{\Cat}{\defcat{Cat}}
\newcommand{\Grph}{\defcat{Grph}}
\newcommand{\Cart}{\defcat{Cart}}

\newcommand{\PolyNat}{\defcat{Poly}_\catE}
\newcommand{\PolyCart}{\defcat{Poly}_\catE^{\mathrm{c}}}
  \newcommand{\DblPolyNat}{\PolyNat}
  \newcommand{\DblPolyCart}{\PolyCart}
\newcommand{\PolyFun}{\defcat{PolyFun}_\catE}
\newcommand{\PolyFunCart}{\defcat{PolyFun}_\catE^{\mathrm{c}}}
  \newcommand{\DblPolyFun}{\PolyFun}
  \newcommand{\DblPolyFunCart}{\PolyFunCart}

\newcommand{\PolyMnd}{\defcat{PolyMnd}}

\newcommand{\PolyEnd}{\defcat{PolyEnd}}
\newcommand{\Span}{\defcat{Span}}
\newcommand{\PSpan}{P\text{-}\Span}
\newcommand{\PPSpan}{P'\text{-}\Span}
\newcommand{\PMultiCat}{P\text{-}\defcat{Multicat}}

\newcommand{\catfont}{\mathscr}
\newcommand{\catC}{\catfont{C}}
\newcommand{\catD}{\catfont{D}}
\newcommand{\catE}{\catfont{E}}
\newcommand{\catL}{\catfont{L}}

\newcommand{\Palg}{P\text{-}\mathrm{alg}}

\newcommand{\slice}[1]{\catE/{#1}}
\newcommand{\Hom}{\mathrm{Hom}}

\newcommand{\Id}[1]{\mathrm{Id}_{#1}}
\newcommand{\id}[1]{\mathrm{Id}_{#1}}
\newcommand{\ladj}[1]{\Sigma_{#1}}
\newcommand{\pbk}[1]{\Delta_{#1}}
\newcommand{\radj}[1]{\Pi_{#1}}

\newcommand{\el}{\operatorname{el}}
\newcommand{\FinSet}{\defcat{FinSet}}
\newcommand{\comma}{\raisebox{1pt}{$\downarrow$}}

\newcommand{\isopil}{\stackrel{\raisebox{0.1ex}[0ex][0ex]{\(\sim\)}}%
{\raisebox{-0.15ex}[0.28ex]{\(\rightarrow\)}}}

\newcommand{\smallprod}[2]{\overset{#2}{\underset{#1}{\scriptstyle{\prod}}}}

\newcommand{\ext}{\operatorname{Ext}}

\newcommand{\poly}[1]{P_{#1}}
\newcommand{\extension}[1]{\poly{#1}}
\newcommand{\exppoly}[1]{P({#1})}
\newcommand{\basechange}[1]{{#1}^{*}}
\newcommand{\cobasechange}[1]{{#1}_{!}}

\author[Gambino]{Nicola Gambino}
 \address{Dipartimento di Matematica e Applicazioni, Universit\`a di Palermo, via Archirafi 34, 90133   Palermo, Italy \& School of Mathematics, The University of Manchester, Oxford Road, Manchester M13   9PL, England }
 \email{nicola.gambino@gmail.com}

 \author[Kock]{Joachim Kock}
 \address{Departament de Matem\`atiques, Universitat Aut\`onoma de Barcelona, 08193 Bellaterra  
 (Barcelona), Spain}
 \email{kock@mat.uab.cat}

\title{Polynomial functors and polynomial monads} 

\date{March 6, 2010}



\begin{document}

\begin{abstract}
  We study polynomial functors over locally cartesian closed categories.
  After setting up the basic theory, we show how polynomial functors
  assemble into a double category, in fact a framed bicategory.  We show
  that the free monad on a polynomial endofunctor is polynomial.  The 
  relationship with operads and other related notions is explored.
\end{abstract}

\maketitle


\section*{Introduction}

\noindent \textbf{Background.} Notions of polynomial functor have proved useful
in many areas of mathematics, ranging from algebra~\cite{MacdonaldI:symfhp,
KontsevichM:defaoo} and topology~\cite{BissonT:dyelab, PirashviliT:polfff} to
mathematical logic~\cite{GirardJ:norfps, MoerdijkI:weltc} and theoretical
computer science~\cite{JayCB:shatp, AbbottM:catc, HancockP:intpwf}.  The present
paper deals with the notion of polynomial functor over locally cartesian closed
categories.  Before outlining our results, let us briefly motivate this
level of abstraction.
   
Among the devices
used to organise and manipulate numbers, polynomials are ubiquitous.  While
formally a polynomial is a sequence of coefficients, it can be viewed also as a
function, and the fact that many operations on polynomial functions, including
composition, can be performed in terms of the coefficients alone is a crucial
feature.  The idea of polynomial functor is to lift the machinery of polynomials
and polynomial functions to the categorical level.
An obvious notion results from letting the category of finite sets take the
place of the semiring of natural numbers, and defining polynomial functors to be
functors obtained by finite combinations of disjoint union and cartesian
product.  It is interesting and fruitful to allow infinite sets.  One reason is
the interplay between inductively defined sets and polynomial functors.  For
example, the set of natural numbers can be characterised as the least solution
to the polynomial equation of sets
\[
X \iso 1 + X \, ,
\]
while the set of finite
planar trees appears as least solution to the equation
\[
X \iso 1+\sum_{n\in\Nat}X^n .
\]
Hence, one arrives at considering as polynomial functors on the
category of sets all the functors of the form
\begin{equation}\label{equ:intro-1var}
X \mapsto \sum_{a\in A}X^{B_{a}} \, ,
\end{equation}
where $A$ is a set and $(B_{a}\mid a\in A)$ is an $A$-indexed family
of sets, which we represent as a map $f: B
\rightarrow A$ with $B_a=f^{-1}(a)$.
It is natural to study also polynomial functors in many variables.  A
$J$-indexed family of polynomial functors in $I$-many variables has the form
\begin{equation} \label{equ:basicpol}
(X_{i}\mid i\in I) \mapsto
\big( \sum_{a\in A_{j}} \prod_{b\in B_{a}} X_{s(b)} \mid j \in J \big) \,,
\end{equation}
where the indexing refers to the diagram of sets
\begin{equation} \label{equ:basicconf}
\xymatrix { I & \ar[l]_s B \ar[r]^f & A \ar[r]^t & J \,.}
\end{equation}
This expression reduces to \eqref{equ:intro-1var} when $I$ and $J$ are singleton
sets.  The functor specified in \eqref{equ:basicpol} is the composite of three
functors: pullback along $s$, the right adjoint to pullback along $f$, and the
left adjoint to pullback along $t$.  The categorical properties of these basic
types of functors allow us to manipulate polynomial functors like
\eqref{equ:basicpol} in terms of their representing
diagrams~\eqref{equ:basicconf}; this is a key feature of the present approach to
polynomial functors.

Although the theory of polynomial functors over $\Set$ is already rich and
interesting, one final abstraction is due: we may as well work 
in any category
with finite limits where pullback functors have both adjoints.  These are the
locally cartesian closed categories, and we develop the theory in this setting,
applicable not only to some current developments in operad theory and
higher-dimensional algebra~\cite{KockJ:polft, KockJ:polfo}, but also in
mathematical logic~\cite{MoerdijkI:weltc}, and in theoretical computer
science~\cite{AbbottM:catc,HancockP:intpwf}.
We hasten to point out that since the category of vector spaces is not locally
cartesian closed, our theory does not immediately apply to various notions of
polynomial functor that have been studied in that context
\cite{MacdonaldI:symfhp, PirashviliT:polfff}.  The precise relationship is under
investigation.

\bigskip

\noindent \textbf{Main results.} Our general goal is to present a mathematically
efficient account of the fundamental properties of polynomial functors over
locally cartesian closed categories, which can serve as a reference for further
developments. 
With this general aim, we begin our exposition by including some known
results that either belong to folklore or were only available in the computer
science literature (but not in their natural generality), giving them a unified treatment and streamlined proofs.
These results mainly concern the diagram representation of strong natural
transformations between polynomial functors, and versions of some 
of these results can be
found in Abbott's thesis~\cite{AbbottM:phd}.
  
Having laid the groundwork, our first main result is to assemble polynomial
functors into a double category, in fact a framed bicategory in the sense of
Shulman~\cite{ShulmanM:frabmf}, hence providing a convenient and precise way of
handling the base change operation for polynomial functors.  There are two
biequivalent versions of this framed bicategory: one is the strict framed
$2$-category of polynomial functors, the other is the (nonstrict) bicategory of
their representing diagrams.
 
Our second main result states that the free monad on a polynomial functor is a
polynomial monad.  This result extends to general polynomial functors the
corresponding result for polynomial functors in a single
variable~\cite{GambinoN:weltdp} and for finitary polynomial functors on the
category of sets~\cite{KockJ:polft,KockJ:polfo}.  We also observe that free 
monads enjoy a double-categorical universal property which is stronger than the 
bicategorical universal property that a priori characterises them.
  
The final section gives some illustration of the usefulness of the
double-category viewpoint in applications.  We give a purely diagrammatic comparison between
Burroni $P$-spans~\cite{BurroniA:tcat}, and polynomials over $P$ (for $P$ a 
polynomial monad).  This
yields in turn a concise equivalence between polynomial monads over $P$ and
$P$-multicategories~\cite{BurroniA:tcat, LeinsterT:higohc},
with base change (multifunctors) conveniently
built into the theory.  Operads are a special case of this.

\bigskip

\noindent \textbf{Related work.} Polynomial functors and closely related notions
have been reinvented several times by workers in different contexts, unaware of
the fact that such notions had already been considered elsewhere.  To help
unifying the disparate developments, we provide many pointers to the literature,
although surveying the different developments in any detail is outside the scope
of this paper.

We should say first of all that our notions of polynomial and polynomial functor
are almost exactly the same as the notions of container and container functor
introduced in theoretical computer science by Abbott, Altenkirch and
Ghani~\cite{AbbottM:phd, AbbottM:catc, Abbott-Altenkirch-Ghani:nested,
Abbott-Altenkirch-Ghani:strictly-positive} to provide semantics for recursive
data types, and studied further in \cite{Morris-Altenkirch:lics09}.  The
differences, mostly stylistic, are explained in Paragraph~\ref{para:containers}.
A predecessor to 
containers were the shapely types of Jay and Cockett~\cite{JayCB:shatp} which we
revisit in Paragraphs~\ref{para:shaf}--\ref{para:shat}.  
The importance of polynomial functors in dependent type theory was first observed by
Moerdijk and Palmgren~\cite{MoerdijkI:weltc}, cf.~Paragraph~\ref{para:W}.
Their polynomial functors are what we call polynomial
functors in one variable.  

The use of polynomial functors in program semantics goes back at least to Manes
and Arbib~\cite{ManesE:algaps}, and was recently explored from a different
viewpoint under the name `interaction systems' in the setting of dependent type
theory by Hancock and Setzer~\cite{HancockP:intpwf} and by
Hyvernat~\cite{HyvernatP:logiis}, where polynomials are also given a
game-theoretic interpretation.  The morphisms there are certain bisimulations,
more general than the strong natural transformations used in the present work.

Within category theory, many related notions have been studied.  In
Paragraph~\ref{para:six} we list six equivalent characterisations of polynomial
functors over $\Set$, and briefly comment on the contexts of the related notions:
familially representable functors of Diers~\cite{DiersY:catl} and
Carboni-Johnstone~\cite{CarboniA:conlfr} (see also \cite[App.~C]{LeinsterT:higohc}), and local right adjoints of
Lamarche~\cite{LamarcheF:modpwc}, Taylor~\cite{TaylorP:quadgl}, and
Weber~\cite{WeberM:genmpr,WeberM:fam2fp}, which notion in the present setting is 
equivalent to the notion of parametric right adjoint of 
Street~\cite{StreetR:pettg}.  We also comment on the relationship
with species and analytic functors \cite{JoyalA:fonaes,BergeronF:comsts}, and
with Girard's normal functors~\cite{GirardJ:norfps}.

Tambara~\cite{TambaraD:mult} studied a notion of polynomial motivated by
representation theory and group cohomology, where the three operations are,
respectively, `restriction', `trace' (additive transfer), and `norm'
(multiplicative transfer).  In Paragraph~\ref{para:tambara}, we give an
algebraic-theory interpretation of one of his discoveries.  Further
study of \myemph{Tambara functors} has been carried out by 
Brun~\cite{BrunM:witvtf}, with applications to Witt vectors.

Most of the results of this paper generalise readily from locally cartesian
closed categories to cartesian closed categories, as we briefly explain in
\ref{generalisation}, if just the `middle maps' $f:B
\to A$ are individually required to be exponentiable.  This generalisation is
useful: for example, covering maps are exponentiable in the category of
compactly generated Hausdorff spaces, and in this way our theory also includes
the notion of polynomial functor used by Bisson and Joyal~\cite{BissonT:dyelab}
to give a geometric construction of Dyer-Lashof operations in bordism.
 
The name polynomial functor is often given to endofunctors of the category of
vector spaces involving actions of the symmetric groups, cf.~Appendix A of
Macdonald's book~\cite{MacdonaldI:symfhp}, a basic ingredient in the
algebraic theory of operads \cite{KontsevichM:defaoo}.  The truncated version of
such functors is a basic notion in functor cohomology, cf.~the survey of
Pirashvili~\cite{PirashviliT:polfff}.  As mentioned, these developments
are not covered by our theory in its present form.

\medskip

This paper was conceived in parallel to~\cite{KockJ:polft,KockJ:polfo}, to take
care of foundational issues.  Both papers rely on the double-categorical
structures described in the present paper, and freely blur the distinction
between polynomials and polynomial functors, as justified in
Section~\ref{sec:morphisms} below.  The paper \cite{KockJ:polfo} uses polynomial
functors to establish the first purely combinatorial characterisation of the
opetopes, the shapes underlying several approaches to higher-dimensional
category theory \cite{LeinsterT:surdnc}, starting with the work of Baez and
Dolan~\cite{BaezJ:higda3}.  In \cite{KockJ:polft}, a new tree formalism based on
polynomial functors is introduced, leading to a nerve theorem charactersising
polynomial monads among presheaves on a category of trees.

\bigskip

\noindent \textbf{Outline of the paper.} In Section~\ref{sec:polf} we recall the
basic facts needed about locally cartesian closed categories, introduce
polynomials and polynomial functors, give basic examples, and show that
polynomial functors are closed under composition.  We also summarise the known
intrinsic characterisations of polynomial functors in the case $\catE=\Set$.  In
Section~\ref{sec:morphisms} we show how strong natural transformations between
polynomial functors admit representation as diagrams connecting the polynomials.
In Section~\ref{sec:double} we assemble polynomial functors into a double
category, in fact a framed bicategory.  In Section~\ref{sec:polmnd} we recall a
few general facts about free monads, and give an explicit construction of the
free monad on a polynomial endofunctor, exhibiting it as a polynomial monad.
Section~\ref{sec:spans-ope} explores, in diagrammatic terms, the relationship
between polynomial monads, multicategories, and operads.

\bigskip

\noindent \textbf{Acknowledgments.} Both authors have had the privilege of being
mentored by Andr\'e Joyal, and have benefited a lot from his generous guidance.
In particular, our view on polynomial functors has been shaped very much by his
ideas, and the results of Section~\ref{sec:morphisms} we essentially learned
from him.  We also thank \mbox{Anders} Kock and Mark Weber for numerous
helpful discussions.  Part of this work was carried out at the CRM in
Barcelona during the special year on Homotopy Theory and Higher Categories; we
are grateful to the CRM for excellent working conditions and for support for the
first-named author.  The second-named author acknowledges support from research
grants MTM2006-11391 and MTM2007-63277 of the Spanish Ministry for Science and
Innovation.

\section{Polynomial functors}
\label{sec:polf}

\begin{para}
  Throughout we work in a locally cartesian closed category $\catE$, assumed to
  have a terminal object~\cite{FreydP:aspt}.  In Section~\ref{sec:polmnd} we 
  shall furthermore assume that $\catE$ has sums and that these are disjoint
  (cf.~also \cite{MoerdijkI:weltc}),
  but we wish to stress that the basic theory of polynomial functors 
  (Section~\ref{sec:polf}, \ref{sec:morphisms}, \ref{sec:double} and 
  \ref{sec:spans-ope}) does not depend on this assumption.
  For $f : B \rightarrow A$
  in~$\catE$, we write $\pbk{f} : \slice{A} \rightarrow \slice{B}$ for pullback
  along $f$.  The left adjoint to $\pbk{f}$ is called the \myemph{dependent sum}
  functor along $f$ and is denoted $\ladj{f} : \slice{B} \rightarrow \slice{A}$.
  The right adjoint to $\pbk{f}$ is called the \myemph{dependent product}
  functor along $f$, and is denoted $\radj{f} : \slice{B} \rightarrow
  \slice{A}$.  We note that
  both unit and counit for the adjunction $\ladj f
  \vdash \pbk f$ are cartesian natural transformations (i.e.~all their 
  naturality squares are cartesian), whereas the unit and counit for $\pbk f \vdash
  \radj f$ are generally not cartesian.

  Following a well-established
  tradition in category theory~\cite{MakkaiM:firocl}, we will use the internal
  logic of $\catE$ to manipulate objects and maps of $\catE$ syntactically
  rather than diagrammatically, when this is convenient. This 
  internal language is the extensional dependent type theory 
  presented in~\cite{SeelyR:loccc}.  
  In the internal language, an object $X\to A$ of $\slice A$ is written as $(X_a 
  \mid a\in A)$, and the three functors 
  associated to $f: B \to A$ take the form
  \begin{align*}
  \pbk f (X_a\mid a\in A) \ & = \ ( X_{f(b)}\mid b\in B)
\\[5pt]
\ladj f (Y_b \mid b\in B) \ & = \ (\sum_{b\in B_a} X_b \mid a\in A)
\\
    \radj f (Y_b \mid b\in B) \ & = \  (\prod_{b\in B_a} X_b \mid a\in A) \,.
\end{align*}
\end{para} 

\begin{para}\label{para:BC-distr}
  We shall make frequent use of the Beck-Chevalley isomorphisms and
  of the distributivity law of dependent sums over dependent
  products~\cite{MoerdijkI:weltc}.  Given a cartesian square
  \[
  \xymatrix {
  \cdot \drpullback \ar[r]^g \ar[d]_u & \cdot \ar[d]^v \\
  \cdot \ar[r]_f & \cdot
  }
  \]
  the Beck-Chevalley isomorphisms are
  \[
    \ladj g \, \pbk u \iso \pbk v \, \ladj f \qquad \text{and} \qquad
    \radj g \, \pbk u \iso \pbk v \, \radj f \,.
    \]

  Given maps $C \stackrel u \longrightarrow B \stackrel f \longrightarrow A$,
  we can construct the diagram
  \begin{equation}\label{distr-diag}
  \xymatrixrowsep{40pt}
  \xymatrixcolsep{27pt}
  \vcenter{\hbox{
  \xymatrix @!=0pt {
  &N \drpullback \ar[rr]^g \ar[ld]_e \ar[dd]^{w=\pbk f(v)}&& M \ar[dd]^{v=\radj f (u)} \\
  C \ar[rd]_u && & \\
  &B \ar[rr]_f && A\,,
  }}}
  \end{equation}
  where $w = \pbk f\, \radj f(u)$ and $e$ is the counit
  of $\pbk{f} \adjoint \radj{f}$.
  For such diagrams the following
  distributive law holds:
  \begin{equation}\label{distr-law}
  \radj f \, \ladj u \iso \ladj v \, \radj g \, \pbk e \,.
  \end{equation}
  In the internal language, the distributive law reads
  \begin{eqnarray}
  \Big( \prod_{b\in B_a} \sum_{c\in C_b} X_c \mid a\in A \Big) &\iso&
  \Big( \sum_{m\in M_a} \prod_{n\in N_m} X_{e(n)} \mid a\in A \Big) \notag \\
 &\iso& \Big(\!\!\!\sum_{m\in\!\!\!\smallprod{b\in B_a}{}\!\!\!C_b} \,\prod_{b\in B_a} X_{m(b)} \mid a\in A \Big) \,.
 \label{equ:distr}
  \end{eqnarray}
\end{para}

\begin{para}
  We recall some basic facts about enrichment,
  tensoring, and strength~\cite{KellyG:bascec,KockA:strfmm}.
    For any object $a:A\to I$ in $\slice I$, the diagram
  $A \stackrel a \to I \stackrel u \to 1$
  defines a pair of adjoint functors
  \[
  \ladj a \pbk a \pbk u \adjoint \radj u \radj a \pbk a \,.
  \]
  The right adjoint provides enrichment of $\slice I$ over $\catE$ by
  setting
  \[
  \underline{\Hom}(a,x) = \radj u \radj a \pbk a (x) \in \catE \,,  \qquad x\in \slice I\,.
  \]
  The left adjoint makes $\slice I$ tensored over $\catE$ by setting
  \begin{equation}\label{tensoring}
  K \otimes a = \ladj a \pbk a \pbk u (K) \in \slice I \,, \qquad K \in \catE \,.
  \end{equation}
  Explicitly, $K \otimes a$ is the object $K \times A \to A \to I$.
  In the internal language, the formulae are (for $a:A\to I$ and 
  $x:X\to I$ in $\slice I$):
  $$
  \underline{\Hom}(a,x) = \prod_{i\in I} X_i^{A_i}\,,
  \qquad\quad 
  K \otimes a = (K \times A_i \mid i\in I)\,.
  $$

  Recall that a tensorial strength~\cite{KockA:strfmm} on a functor
  $F:\catD\to\catC$ between categories tensored over $\catE$ is a family of maps
  $$
  \tau_{K,a} : K \otimes F(a) \to F(K\otimes a)
  $$
  natural in $K\in \catE$ and in $a\in \catD$, and satisfying two axioms
  expressing an associativity and a unit condition.  A
  natural transformation between strong functors is called strong if it is
  compatible with the given strengths.  When $\catE$ is cartesian closed, giving
  a tensorial strength is equivalent to giving an enrichment, and a natural
  transformation is strong if and only if it is enriched.
  
  For any $f: B \to A$, there is a canonical strength on
  each of the three functors $\pbk f$, $\ladj f$, and $\radj f$:
  writing out using Formula~\eqref{tensoring} it is easily seen that
  the strength on $\pbk f$ is essentially a Beck-Chevalley isomorphism,
  the strength of $\ladj f$ is essentially trivial, whereas the strength of 
  $\radj f$ depends on distributivity and is essentially an instance of the
  unit for the $\pbk{} \adjoint \radj{}$ adjointness.
  It is also direct to verify that
  the natural transformations given by the units and counits for the
  adjunctions, as well as those expressing pseudo-functoriality of pullback and 
  its adjoints, are all strong
  natural transformations.  We shall work with strong functors and strong
  natural transformations, as a convenient alternative to the purely enriched
  viewpoint.
\end{para}

\begin{para} 
    We define a \myemph{polynomial} over $\catE$ to be a diagram $F$ in $\catE$
    of shape
\begin{equation}
\label{equ:poly}
\xymatrix{
I  & B \ar[l]_s \ar[r]^f & A \ar[r]^t & J \, . }
\end{equation}
We define $\poly{F}:\slice I \to \slice J$ as the composite
\[
\xymatrix{
\slice{I} \ar[r]^{\pbk{s}} & \slice{B} \ar[r]^{\radj{f}} & \slice{A} \ar[r]^{\ladj{t}} & \slice{J} \, . }
\]
We refer to $\poly{F}$ as the \myemph{polynomial functor} associated to $F$, or
the \myemph{extension} of $F$, and say that $F$ \myemph{represents} $\poly F$.
In the internal language of $\catE$, the functor $\poly{F}$ has the expression
\[
\poly{F}(X_i \mid  i \in I) =
\Big( \sum_{a \in A_j } \prod_{b \in B_a  } X_{s(b)} \mid j \in J \Big) \, .
\]
By a \myemph{polynomial functor} we understand any functor
isomorphic to the extension of a polynomial.  
The distinction between polynomial and polynomial functor is similar to the
usage in elementary algebra, where a polynomial defines a polynomial function.
The bare polynomial is an abstract configuration of exponents and coefficients
which can be interpreted by extension as a function.  This extension is of
course a crucial aspect of polynomials, and conversely it is a key feature of
polynomial functions that they can be manipulated in terms of the combinatorial
data.  A similar interplay characterises the theory of polynomial functors.  We
shall shortly establish a result justifying the blur between polynomials and
polynomial functors; only in the present paper do we insist on the distinction.
\end{para}

\begin{para}
When $I = J = 1$,  a polynomial is essentially given by a single map $B 
\rightarrow A$, and the extension reduces to
\[
P(X) = \sum_{a \in A} X^{B_a} \, .
\]
Endofunctors of this form, simply called polynomial functors in
\cite{MoerdijkI:weltc}, will be referred to here as \myemph{polynomial functors in a
single variable}.
\end{para}

\begin{examples}\label{ex:id} \, \hfill
\begin{enumerate}[(i)]
 \item The identity functor $\Id{} : \slice{I} \to \slice{I}$ is polynomial,
it is represented by $$I \stackrel = \leftarrow  I \stackrel = 
\rightarrow I \stackrel = 
\rightarrow I\,.$$
 \item  If $\catE$ has an initial object $\emptyset$, then for any $A\in \slice J$, 
 the constant functor
 $\slice{I} \to \slice{J}$ with value $A$ is polynomial, represented by
 $$
 I\stackrel s \leftarrow \emptyset\rightarrow A\rightarrow J\,.
 $$
 (Indeed already $\pbk s$ is constant $\emptyset$.)
\end{enumerate}
\end{examples}

\begin{example}\label{ex:linear}
  A span $I \stackrel s \leftarrow M \stackrel t \rightarrow J$ can be regarded as
  a polynomial
  $$
  I \stackrel s \leftarrow M \stackrel = \to M  \stackrel t \rightarrow J .
  $$
  The associated polynomial functor
  $$
  \poly{M}(X_{i}\mid i\in I) = \big(\sum_{m\in M_{j}}X_{s(m)}\mid j\in J\big)
  $$
  is called a \myemph{linear} functor, since it is the formula for matrix 
  multiplication, and since $\poly{M}$ preserves sums.
  Hence polynomials can be seen as a natural `non-linear'
  generalisation of spans.
\end{example}

\begin{example}
  Let $C = (C_0 \overset s {\underset {t}{\leftleftarrows}} C_1)$ be a category object in $\catE$.
  The polynomial
  $$
  C_0 \stackrel s \leftarrow C_1 \stackrel = \to C_1 \stackrel t \to C_0
  $$
  represents the polynomial (in fact linear) endofunctor 
  $\slice{C_0} \to \slice{C_0}$
  which gives the free internal
  presheaf on a $C_0$-indexed family~\cite[\S V.7]{MacLaneS:shegl}.
\end{example}

\begin{example}\label{ex:freemonoid}
The free-monoid monad, also known as the word monad or the list monad,
  \begin{eqnarray*}
    M : \Set & \longrightarrow & \ \Set  \\
    X & \longmapsto & \sum_{n\in\Nat} X^n
  \end{eqnarray*}
  is polynomial, being represented by the diagram
 \[
 \xymatrix{ 1  & \ar[l]   \Nat'     \ar[r]  & \Nat \ar[r] &  1\,,} 
 \] 
 where $\Nat'\to \Nat$ is such that the fibre over $n$ has cardinality $n$, as given for example
 by the second projection from  $\Nat' = \{ (i,n) \in \Nat\times\Nat \mid i<n \}$.
\end{example}

\begin{example}\label{ex:trees}
  (Cf.~\cite{KockJ:polft}.)  A rooted tree defines a polynomial in $\Set$:
 \[
 \xymatrix{
 A & M \ar[l]_{s} \ar[r]^f & N \ar[r]^t & A }
 \]
  where $A$ is the set of edges, $N$ is the set of nodes, and $M$ is the set of
  nodes with a marked incoming edge.  The map $t$ returns the outgoing edge of
  the node, the map $f$ forgets the marked edge, and the map $s$ returns the
  marked edge.  It is shown in \cite{KockJ:polft} that every polynomial is a
  colimit of trees in a precise sense.
\end{example}

\begin{para} \label{para:comp}
  We now define the operation of substitution of polynomials, and show that the
  extension of substitution is composition of polynomial functors, as expected.
  In particular, the composite of two polynomial functors is again polynomial.
  Given polynomials
\[
\xymatrix{
 &  B \ar[r]^f  \ar[dl]_{s} & A \ar[dr]^{t} & \\
 I \ar@{}[rrr]|{F} & & & J } \qquad
\xymatrix{
 & D \ar[r]^g \ar[dl]_{u} & C \ar[dr]^{v}   \\
J \ar@{}[rrr]|{G}& & & K }
\]
we say that $F$ is a polynomial \myemph{from} $I$ \myemph{to} $J$ (and $G$ from 
$J$ to $K$), and
we define $G\circ F$, the \myemph{substitution} of $F$ into $G$, to be the polynomial
$I \leftarrow N \to M \to K$ constructed via this diagram:
\begin{equation}
\label{equ:compspan}
\xycenter{
  &  &  & N \ar[dl]_{n} \ar[rr]^{p} \ar@{}[dr] |{(iv)}
&  & D' 
\ar[dl]^{\varepsilon} \ar[r]^{q} \ar@{}[ddr] |{(ii)} &  
M  \ar[dd]^{w} &  \\
  &  & B' \ar[dl]_{m} \ar[rr]^{r} \ar@{} [dr]|{(iii)}
&  & A'  \ar[dr]^{k} \ar[dl]_{h} 
\ar@{} [dd] |{(i)}
&  &  &  \\
   & B \ar[rr]^f \ar[dl]_{s} & & A \ar[dr]_{t} &   & D \ar[dl]^{u} 
\ar[r]^{g} & C \ar[dr]^{v} &    \\
I  &    & &   & J &   &   & K   }
\end{equation} 
Square $(i)$ is cartesian, and $(ii)$ is a distributivity diagram
like \eqref{distr-diag}: $w$ is obtained
by applying $\radj{g}$ to $k$, and $D'$ is the pullback of $M$ along $g$.
The arrow $\varepsilon: D' \to A'$ is the $k$-component of the counit of
the adjunction $\ladj g \adjoint \pbk g$.
Finally, the squares $(iii)$ and $(iv)$ are
cartesian.
\end{para}

\begin{proposition}\label{thm:subst}
  There is a natural isomorphism
  $$
  \poly{G\circ F} \iso \poly{G} \circ \poly{F} .
  $$
\end{proposition}

\begin{proof}
  Referring to Diagram~\eqref{equ:compspan} we have 
  the following chain of natural isomorphisms:
\begin{eqnarray*}
\extension{G} \circ \extension{F} & = & \ladj{v} \, \radj{g} \, \pbk{u}  \; \ladj{t}  \, \radj{f}   
\, \pbk{s} \\
& \iso & \ladj{v} \, \radj{g}  \, \ladj{k}  \, \pbk{h}  \, \radj{f}  \, \pbk{s}\\
 & \iso & 
\ladj{v} \, \ladj{w}  \, \radj{q}  \, \pbk{\varepsilon}  \, \pbk{h}  \, \radj{f} \, \pbk{s}\\
& \iso &
\ladj{v} \, \ladj{w}  \, \radj{q}  \, \radj{p}  \, \pbk{n}  \, \pbk{m} \, \pbk{s}\\
& \iso & 
\ladj{(v\, w)} \, \radj{(q \, p)}  \, \pbk{(s\, m\, n)}\, \\
& = & \extension{G \circ F}\,.
\end{eqnarray*}
Here we used the Beck-Chevalley isomorphism for the cartesian square
$(i)$, the distributivity law for $(ii)$, Beck-Chevalley isomorphism
for the cartesian squares $(iii)$ and $(iv)$, and finally pseudo-functoriality 
of the pullback functors and their adjoints. 
\end{proof}

\begin{para} 
  Let us also spell out the composition in terms of the 
internal language, to highlight the
substitutional aspect. By definition, the composite functor is given by
\[
\extension{G} \circ \extension{F}  (X_i \mid i \in I) = 
\Big( \sum_{c \in C_k} \prod_{d \in D_c} \sum_{a \in A_{u(d)} } 
\prod_{b \in B_a} X_{s(b)} \mid k \in K \Big)
\, . 
\]
For fixed $c\in C$, 
by distributivity~\eqref{equ:distr}, we have
\[
 \prod_{d \in D_c} \sum_{a \in A_{u(d)} } \prod_{b \in B_a} X_{s(b)}  
\iso 
 \sum_{m \in M_c} \prod_{d \in D_c} \prod_{b \in B_{m(d)} } X_{s(b)} \,,
\]
where we have put
\[
M_c =  \prod_{d \in D_c} A_{u(d)} \,,
\]
the $w$-fibre over $c$ in Diagram~\eqref{equ:compspan}.
If we also put, for $m\in M_c$,
\[
N_{(c,m)} =  \sum_{d \in D_c} B_{m(d)} \,, 
\]
the $(q\circ p)$-fibre over $m\in M_c$,
we can write
\[
\sum_{m \in M_c} \prod_{d \in D_c} \prod_{b \in B_{m(d)} } X_{s(b)}  \iso
\sum_{m \in M_c} \prod_{(d,b) \in N_{(c,m)} } X_{s(b)} \,.
\]
Summing now over $c\in C_k$, for $k\in K$, we conclude
\[
\extension{G} \circ \extension{F} (X_i \mid i \in I)  \iso \Big(
\sum_{(c,m) \in M_k} \prod_{(d,b) \in N_{(c,m)} } X_{s(b)} \mid k \in K 
\Big) \, ,
\]
(where $M_k = \sum_{c \in C_k} M_c$ is the $(v\circ w)$-fibre over $k\in K$).
\end{para}

\begin{corollary} 
  The class of polynomial functors is the smallest class of functors between
  slices of $\catE$ containing the pullback functors and their adjoints,
  and closed under composition and natural isomorphism. \qed
\end{corollary} 

\begin{proposition} 
  Polynomial functors have a natural strength.
\end{proposition} 

\begin{proof}
  Pullback functors and their adjoints have a canonical strength.
\end{proof}

\begin{proposition} Polynomial functors preserve connected limits.
    In particular, they are cartesian.
\end{proposition} 

\begin{proof} Given a diagram as in~\eqref{equ:poly}, the functors $\pbk{s} :
\slice{I} \rightarrow \slice{B}$ and $\radj{f} : \slice{B} \rightarrow \slice{A}$
preserve all limits since they are right adjoints. A direct calculation shows that
also the functor $\ladj{t} : \slice{A}
\rightarrow \slice{J}$ preserves connected limits 
\cite{CarboniA:conlfr}. 
\end{proof} 

\begin{para}\label{generalisation}
  In this paper we have chosen to work with locally cartesian closed categories,
  since it is the most natural generality for the theory.  However, large parts
  of the theory make sense also over cartesian closed categories, by considering
  only polynomials for which the `middle map' $f: B \to A$ is exponentiable, or
  belongs to a subclass of the exponentiable maps having the same stability
  properties to ensure that Beck-Chevalley, distributivity, and composition of
  polynomial functors work just as in the locally cartesian closed case.
  Further results about polynomial functors in this generality can be deduced
  from the locally cartesian closed theory by way of the Yoneda embedding $y:
  \catE \to \widehat \catE$, where $\widehat \catE$ denotes the category of
  presheaves on $\catE$ with values in a category of sets so big that $\catE$ is
  small relatively to it.  The Yoneda embedding is compatible with slicing and
  preserves the three basic operations, so that basic results about polynomial
  functors in $\catE$ can be proved by reasoning in $\widehat\catE$.  A
  significant example of this situation is the cartesian closed category of
  compactly generated Hausdorff spaces, where for example the covering maps
  constitute a stable class of exponentiable maps.  Polynomial functors in this
  setting were used by Bisson and Joyal~\cite{BissonT:dyelab} to give a
  geometric construction of Dyer-Lashof operations in bordism.  Another example
  is the category of small categories, where the Conduch\'e fibrations are the
  exponentiable maps.  In this setting, an example of a polynomial functor is
  the family functor, associating to a category $X$ the category of families of
  objects in $X$.
\end{para}
  
\begin{para}\label{para:six}
  For the remainder of this section, with the aim of putting the theory of
  polynomial functors in historical perspective, we digress into the special case
  $\catE=\Set$, then make some remarks on finitary polynomial functors, and end
  with finite polynomials.  This material is not needed in the subsequent
  sections.
  
  The case $\catE=\Set$ is somewhat special due to the equivalence $\Set/I
  \simeq \Set^I$, which allows for various equivalent characterisations of
  polynomial functors over $\Set$.
  
  For a functor $P :
  \Set/I \rightarrow \Set/J$, the following conditions are equivalent.
  \begin{enumerate}[(i)]
  \item \label{item:poly} $P$ is polynomial.
  \item \label{item:connected} $P$ preserves connected limits 
  (or, equivalently, pullbacks and cofiltered limits, 
  or equivalently, wide pullbacks).
  \item \label{item:famrepr} $P$ is familially representable (i.e.~a sum of 
  representables).
  \item \label{item:topos} The comma category $(\Set/J) \comma P$ is a presheaf
  topos.
  \item \label{item:lra} $P$ is a local right adjoint (i.e.~the slices of $P$ are right adjoints).
  \item \label{item:generic} $P$ admits strict generic factorisations 
  \cite{WeberM:genmpr}.
  \item \label{item:normalform} Every slice of $\operatorname{el}(P)$ has an 
  initial object (Girard's normal-form property).
  \end{enumerate}
  
  The equivalences $\textrm{(\ref{item:connected})} \Leftrightarrow
  \textrm{(\ref{item:lra})} \Leftrightarrow \textrm{(\ref{item:generic})}$ go
  back to Lamarche~\cite{LamarcheF:modpwc} and Taylor~\cite{TaylorP:quadgl}, who
  were motivated by the work of Girard~\cite{GirardJ:norfps}, cf.~below.  They
  arrived at condition (\ref{item:generic}) as the proper generalisation of
  (\ref{item:normalform}), itself a categorical reformulation of Girard's
  normal-form condition \cite{GirardJ:norfps}.  Below we give a direct proof of
  $\textrm{(\ref{item:poly})} \Leftrightarrow \textrm{(\ref{item:normalform})}$,
  to illuminate the relation with Girard's normal functors.  The
  equivalence~$\textrm{(\ref{item:connected})} \Leftrightarrow
  \textrm{(\ref{item:famrepr})}$ is due to Diers~\cite{DiersY:catl}, and was
  clarified further by Carboni and Johnstone~\cite{CarboniA:conlfr}, who
  established in particular 
  the equivalence~$\textrm{(\ref{item:connected})} \Leftrightarrow
  \textrm{(\ref{item:topos})}$ as part of their treatment of Artin gluing.
  The equivalence~$\textrm{(\ref{item:poly})} \Leftrightarrow
  \textrm{(\ref{item:famrepr})}$ is also implicit in their work, the
  one-variable case explicit.  The equivalence~$\textrm{(\ref{item:poly})}
  \Leftrightarrow \textrm{(\ref{item:lra})}$ was observed by
  Weber~\cite{WeberM:fam2fp}, who also notes that on general presheaf toposes,
  local right adjoints need not be polynomial: for example the free-category
  monad on the category of directed graphs is a local right adjoint but not a
  polynomial functor.
\end{para}

\begin{para}
  A polynomial functor $P:\Set/I \to \Set/J$ is \myemph{finitary} if it
  preserves filtered colimits.  If $P$ is represented by $I \leftarrow B \to A
  \to J$, this condition is equivalent to the map $B \to A$ having finite
  fibres.
\end{para}

\begin{para}
  Recall~\cite{Joyal:1981,BergeronF:comsts} that a \myemph{species}
  is a functor $F: \FinSet_{\text{bij}} \to \Set$,
  or equivalently, a sequence $(F[n]\mid n\in\Nat)$
  of $\Set$-representations of the symmetric groups.
  To a species is associated an \myemph{analytic functor}
  \begin{eqnarray*}
    \Set & \longrightarrow & \Set \\
    X & \longmapsto & \sum_{n\in\Nat} F[n] \times_{\mathfrak S_n} X^n\,.
  \end{eqnarray*}
  Species and analytic functors were introduced by Joyal~\cite{JoyalA:fonaes},
  who also characterised analytic
  functors as those preserving weak pullbacks, cofiltered limits, and filtered colimits.
  It is the presence of group actions that makes the preservation of pullbacks 
  weak, in contrast to the polynomial functors, cf.~(\ref{item:connected}) 
  above.  Species for which the group actions are free are called \myemph{flat}
  species \cite{BergeronF:comsts}; they encode rigid 
  combinatorial structures, and correspond to ordinary generating functions 
  rather than exponential ones.  The analytic functor associated to a flat species preserves
  pullbacks strictly and is therefore the same thing as a finitary polynomial
  functor on $\Set$. Explicitly,
  given a one-variable finitary polynomial functor $P(X)=
  \sum_{a\in A}X^{B_{a}}$ represented by $B \to A$, we can `collect terms': let $A_n$ denote
  the set of fibres of cardinality $n$, then there is a bijection
  $$
  \sum_{a\in A}X^{B_{a}} \iso \sum_{n\in\Nat} A_n \times X^n .
  $$
  The involved bijections $B_{a}\iso n$ are not canonical: the degree-$n$ part 
  of $P$ is rather a $\mathfrak S_n$-torsor, denoted $P[n]$, and we
  can write instead
  \begin{equation}\label{equ:species}
  P(X) \iso \sum_{n\in\Nat} P[n] \times_{\mathfrak S_n} X^n ,
  \end{equation}
  which is the analytic expression of $P$.
  
  As an example of the polynomial encoding of a flat species, consider the
  species $C$ of binary planar rooted trees.  The associated analytic functor
  is
  $$
  X \mapsto \sum_{n\in \Nat} C[n] \times_{\mathfrak S_n} X^n \,,
  $$
  where $C[n]$ is the set of ways to organise an $n$-element set as the set of
  nodes of a binary planar rooted tree; $C[n]$ has cardinality $n!  \, c_{n}$,
  where $c_n$ are the Catalan numbers $1,1,2,5,14,\ldots$
  The polynomial representation is
  $$
  1 \longleftarrow B \longrightarrow A \longrightarrow 1
  $$
  where $A$ is the set of isomorphism classes of binary planar rooted trees,
  and $B$ is the set of isomorphism classes of binary planar rooted trees with
  a marked node. 
\end{para}

\begin{para}
  Girard~\cite{GirardJ:norfps}, aiming at constructing models for lambda
  calculus, introduced the notion of \myemph{normal functor}: it is
  a functor $\Set^I \to \Set^J$ which preserves pullbacks, cofiltered
  limits and filtered colimits, i.e.~a finitary
  polynomial functor.  Girard's interest was a certain normal-form property
  (reminiscent of Cantor's normal form for ordinals), which in modern
  language is (\ref{item:normalform}) above: the normal forms of the functor
  are the initial objects of the slices of its category of elements.
  Girard, independently of \cite{JoyalA:fonaes}, also
  proved that these functors admit a power series expansion, which is just
  the associated (flat) analytic functor.  From Girard's proof we can extract in
  fact a direct equivalence between (\ref{item:poly}) and (\ref{item:normalform})
  (independent of the finiteness
  condition).
  The proof shows that, in a sense, the polynomial representation
  \myemph{is} the normal form.
  For simplicity we treat only the one-variable case.
\end{para}

\begin{proposition}\label{thm:normalform}
  A functor $P:\Set\to\Set$ is polynomial if and only if every slice of
  $\el(P)$ has an initial object.
\end{proposition}

\begin{proof}
  Suppose $P$ is polynomial, represented by $B \to A$.  
  An element of $P$ is a triple
  $(X,a,s)$, where $X$ is a set, $a\in A$, and $s:B_a \to X$.  
  The set of connected components of $\el(P)$ is in bijection with the set
  $P(1) = A$.  For each element
  $a\in A=P(1)$, it is clear that
  the triple $(B_a, a, \id{B_a})$ is an initial object of the slice
  $\el(P)/(1,a,u)$, where $u$ is the map to the terminal object.
  These initial objects induce initial objects in 
  all the slices, since every element $(X,a,s)$ has a unique map to $(1,a,u)$.
  
  Conversely, suppose every slice of $\el(P)$ has an initial object; again we
  only need the initial objects of the special slices
  $\el(P)/(1,a,u)$, for $a\in P(1)$.  Put $A = P(1)$.  It
  remains to construct $B$ over $A$ and show that the resulting polynomial
  functor is isomorphic to $P$.  Denote by $(B_a,b)$ the initial object of
  $\el(P)/(1,a,u)$.  Let now $X$ be any set.  The unique map $X \to 1$ induces
  $P(X)\to P(1)=A$, and we denote by $P(X)_a$ the preimage of $a$.  For each
  element $x\in P(X)_a$, the pair $(X,x)$ is therefore an object of the slice
  $\el(P)/(1,a,u)$, so by initiality we get a map $B_a\to X$.  Conversely, given any map
  $\alpha:B_a \to X$, define $x$ to be the image under $P(\alpha)$ of the element
  $b$; clearly $x\in P(X)_a$.  These two constructions are easily checked to be
  inverse to each other, establishing a bijection $P(X)_a \iso X^{B_a}$.
  These bijections are clearly natural in $X$, and since
  $P(X)=\sum_{a\in A} P(X)_a$ we conclude that $P$ is isomorphic to the
  polynomial functor represented by the projection map $\sum_{a\in A} B^a \to A$.
\end{proof}

\bigskip

\begin{para}\label{para:tambara}
  Call a polynomial over $\Set$
  \begin{equation}\label{equ:polytamb}
    I \leftarrow B \to A \to J
  \end{equation}
  \myemph{finite} if the four
  involved sets are finite.  Clearly the composite of two finite
  polynomials is again finite.  The category $\mathbb{T}$ whose objects are
  finite sets and whose morphisms are the finite polynomials (up to isomorphism)
  was studied by Tambara~\cite{TambaraD:mult}, 
  in fact in the more general context of finite $G$-sets, for $G$ a finite group. His paper is probably the first to display
  and give emphasis to diagrams like \eqref{equ:polytamb}. Tambara was
  motivated by representation
  theory and group cohomology, where the three operations $\pbk{}$, $\ladj{}$,
  $\radj{}$ are, respectively, `restriction', `trace' (additive transfer), and
  `norm' (multiplicative transfer).  We shall not go into the $G$-invariant
  achievements of 
  \cite{TambaraD:mult}, but
  wish to point out that the following fundamental result about polynomial 
  functors is implicit in Tambara's paper and should be attributed to him.

\begin{theorem}
  The skeleton of $\mathbb{T}$ is the Lawvere theory for commutative
  semi\-rings.
\end{theorem}
\noindent
The point is firstly that $m+n$ is the product of $m$ and $n$ in $\mathbb T$
(this is most easily seen by extension, where it amounts to $\Set/(m+n)\simeq
\Set/m\times\Set/n$).  And secondly that
for the two $\Set$-maps
\[
0 \stackrel e \longrightarrow 1 \stackrel m \longleftarrow 2
\]
the polynomial functor $\ladj m$, considered as a map in $\mathbb T$,
represents addition, $\radj m$ represents multiplication, and $\ladj e$
and $\radj e$ represent the additive and multiplicative neutral elements,
respectively.  Pullback provides the projection for the product in $\mathbb T$,
and is also needed to account for distributivity,
which in syntactic terms involves duplicating elements.
It is a beautiful exercise to use the abstract distributive law \eqref{distr-law}
to compute
$$
\radj m \circ \ladj k
$$
where $k:3\to 2$ is the map pictured as
\raisebox{-1.5pt}{
\begin{texdraw}
  \move (0 0) \fcir f:0 r:1 
  \lvec (8 2) \fcir f:0 r:1
  \move (0 4) \fcir f:0 r:1
  \lvec (8 6) \fcir f:0 r:1 
  \move (0 8) \fcir f:0 r:1 
  \lvec (8 6)
\end{texdraw}
} ,
recovering the distributive law $a(x+y)=ax+ay$ of elementary algebra.
\end{para}

\section{Morphisms of polynomial functors} 
\label{sec:morphisms}

Since polynomial functors have a canonical strength, the natural notion of
morphism between polynomial functors is that of strong natural transformation.
We shall see that strong natural transformations between polynomial functors are
uniquely represented by certain diagrams connecting the polynomials.

\begin{para} \label{para:cart}
    Given a diagram 
\begin{equation}
\label{equ:cartmor}
\xycenter{
F': &I \ar@{=}[d] & B' \drpullback \ar[r]^{f'} \ar[l]_{s'}  
\ar[d]_{\beta} & A' \ar[r]^{t'} \ar[d]_{\alpha} & J \ar@{=}[d] &\\
F: &I & B \ar[r]_{f} \ar[l]^{s} & A \ar[r]_{t} & J &
}
\end{equation}
we define a cartesian strong natural transformation
$\phi : \poly{F'}  \Rightarrow \poly{F}$ by the pasting
diagram
\[
\xymatrix{
\slice{I'} \ar[rr]^{\pbk{s'}}  \ar[dr]_{\pbk{s}} &  \ar@{}[d]|{\iso} & \slice{B'} \ar[rr]^{\radj{f'}} 
\ar@{}[dr]|{\iso} & & \slice{A'} \ar@{}[d]|{\ \Downarrow \, \varepsilon}
\ar[rr]^{\ladj{t'}}  \ar[dr]^{\ladj{\alpha}} &  &  \slice{J'} \\
& \slice{B} \ar[rr]_{\radj{f}}  \ar[ur]_{\pbk{\beta}} & &   \slice{A} \ar@{=}[rr]  \ar[ur]^{\pbk{\alpha}} 
& &  \slice{A} \ar[ur]_{\ladj{t}} & }
\]
It is cartesian and strong
since its constituents are so.

In the internal language of $\catE$, the component of 
$\phi : \poly{F'}  \Rightarrow \poly{F}$ at $X = (X_i \mid i \in I)$ 
is the function
\[
\phi_X :
\Big( \sum_{a' \in A'_j} \prod_{b' \in B'_{a'}} X_{s'(b')} \mid j \in J \Big)    \rightarrow
\Big( \sum_{a \in A_j} \prod_{b \in B_{a}} X_{s(b)} \mid j \in J \Big)
\]
defined by 
\[
\phi_X(a',x' ) =  \big( \alpha(a'), x' \cdot \beta_{a'}^{-1} \big) \, ,
\]
where  $\beta_{a'} : B'_{a'} \rightarrow B_{\alpha(a')}$ is the 
isomorphism determined by the cartesian square in~\eqref{equ:cartmor}. 
\end{para} 

\begin{lemma} \label{thm:carttopoly}
    Let $P:\slice I \to \slice J$ be a polynomial functor.
    If $Q \Rightarrow P$ is a cartesian natural transformation,
    then $Q$ is also a polynomial functor.
\end{lemma}

\begin{proof}
    Assume $P$ is represented by $I \leftarrow B \to A \to J$.
    Construct the diagram
\[
\xycenter{
I   \ar@{=}[d]    & B' \drpullback \ar[r]^{f'}  \ar[l]_{s'}  
\ar[d]_{\beta}  & A'  \ar[r]^{t'}  \ar[d]_{\alpha} & J \ar@{=}[d] \\
I                         & B  \ar[r]_{f}   
\ar[l]^{s}                             & A  \ar[r]_{t} & J   }
\]
by setting $A' = Q(1)$, and taking $\alpha : A' \rightarrow A$ to
be the map $\phi_1 : Q(1) \rightarrow P(1)$, and letting $B'$ 
be the pullback.  The top row represents a polynomial functor $P'$, and the
diagram defines  a
cartesian natural transformation to $P$.
Since $P'$ and $Q$ both have a cartesian natural transformation to $P$ which agree
on the terminal object, they are naturally isomorphic.  Hence $Q$ is polynomial.
\end{proof} 

\begin{para} 
Recall that, for a category $\catC$  with a terminal object $1$
and a category~$\catD$ with pullbacks, the functor 
\begin{eqnarray*}
  {}[\catC,\catD] & \longrightarrow & \catD  \\
  P & \longmapsto & P(1) 
\end{eqnarray*}
is a Grothendieck fibration. The cartesian arrows for this fibration are precisely 
the cartesian natural transformations, while the vertical arrows are the natural 
transformations whose component on $1$ is an identity map.
We refer to such natural transformations as \myemph{vertical natural transformations}.

If $\catC$ and $\catD$ are enriched and tensored, then the above remark carries
over to the case where $[\catC,\catD]$ denotes the category of
strong functors and strong natural transformations.  The verification of this 
involves observing that the cartesian lift of a strong functor has a canonical
strength.
\end{para}

\begin{proposition} \label{thm:grothfib}
Let $I, J \in \catE$. The restriction of the Grothendieck fibration
$[ \slice{I}, \slice{J} ] \rightarrow \slice{J}$ 
to the category of polynomial functors and strong natural transformations 
is  again a Grothendieck fibration.
\end{proposition}

\begin{proof} 
  Lemma~\ref{thm:carttopoly} implies that the cartesian lift of a polynomial
  functor is again polynomial.
\end{proof}

\begin{para} 
  Proposition~\ref{thm:grothfib} implies that every strong natural
  transformation between polynomial functors factors in an essentially unique
  way as a vertical strong natural transformation followed by a cartesian one.  We
  proceed to establish representations of the two classes of strong natural
  transformations between polynomial functors.  The key ingredient is the
  following version of the enriched Yoneda lemma.
\end{para}

\begin{lemma}
  Let $u:I\to 1$ denote the unique arrow in $\catE$ to the terminal object.  For
  any $s: B \to I$ and $s': B' \to I$ in $\slice I$, the natural map
  \[
  \Hom_{\slice I}(s,s') \longrightarrow 
  \defcat{StrNat} (\radj{u} \radj{s'} \pbk{s'},\radj{u} \radj{s} \pbk{s}) 
  \]
  sending an $I$-map $w:B\to B'$ to the composite $\radj{u} \radj{s'} \pbk{s'}
  \stackrel \eta \Rightarrow \radj{u} \radj{s'} \radj{w}\pbk{w} \pbk{s'}
  \iso \radj{u} \radj{s} \pbk{s}$ is a bijection.
\end{lemma}

\begin{proof}
  Just note that $\radj{u} \radj{s} \pbk{s} = \underline{\Hom}_{\slice I}(s, - ) : \slice I
  \to \catE$, and the result is the usual enriched Yoneda lemma
  \cite{KellyG:bascec}, remembering that since $\slice I$ is tensored over
  $\catE$, a natural transformation (between strong functors) is enriched if and
  only if it is strong.
\end{proof}

\begin{para} \label{para:vert}
Given a  diagram
\begin{equation}
\label{equ:repmor}
\xycenter{
F': &I \ar@{=}[d] & B'  \ar[r]^{f'}  \ar[l]_{s'} & A \ar[r]^{t} & J \ar@{=}[d] & \\
F: & I  & B  \ar[r]_{f}  \ar[u]_{w}    \ar[l]^{s}    & A  \ar@{=}[u] \ar[r]_t & 
J & }
\end{equation}
we define a strong natural transformation $\phi : \poly{F'} \Rightarrow \poly{F}$
 by the pasting diagram
\[
\xymatrix{
   & \slice{B'} \ar@{=}[rr] \ar[dr]^{\pbk{w}}  \ar@{}[d]|{\iso} & 
   \ar@{}[d]|{\ \Downarrow \eta}       
    & \slice{B'}  \ar@{}[d]|{\iso}  \ar[dr]^{\radj{f'}}  &                  &                  \\
   \slice{I} \ar[ur]^{\pbk{s'}}          \ar[rr]_{\pbk{s}}          &                            & \slice{B} 
   \ar[ur]^{\radj{w}} \ar[rr]_{\radj{f}} &                    & \slice{A} \ar[r]_{\ladj{t}} & \slice{J}  }
\]
In the internal language,  the  component of $\phi$ at $X = 
(X_i \mid i \in I)$ is given
by the function
\[ 
   \phi_X : \Big( \sum_{a \in A_j} \prod_{b' \in B'_a} X_{u(b)} \mid j \in J \Big) \rightarrow 
   \Big( \sum_{a \in A_j} \prod_{b \in B_a} X_{s(b)} \mid j \in J \Big) 
 \]
 defined by
 \[
 \phi_X (  a, x )  =  \big( a,  x \cdot w_a)
 \]
 Clearly $\phi_1 = \id{A}$, so $\phi$ is vertical for the Grothendieck 
 fibration.
\end{para} 

\begin{proposition} \label{thm:vertmor}
  For $F$ and $F'$ as above, every vertical strong natural transformation
  $\phi: \poly{F'} \Rightarrow \poly{F}$ is uniquely represented by a diagram 
  like \eqref{equ:repmor}.
\end{proposition} 

\begin{proof}
  We already have the outline of the diagram \eqref{equ:repmor}, it remains to
  construct the map $w : B \rightarrow B'$ commuting with the rest.
  Since $w$ must be an $A$-map, we can construct it fibrewise, so we need
  for each $a\in A$ a map $B'_a \to B_a$.  This allows reduction
  to the case $A=1$, and the result is a direct consequence of the above
  Yoneda lemma.
\end{proof}
 
\begin{proposition}\label{thm:cartmor}
  Let $I, J \in \catE$.  Let $F : I \rightarrow J$ and $F' : I \rightarrow J$ be
  polynomials.  Every cartesian strong natural transformation $\phi : \poly{F'}
  \Rightarrow \poly{F}$ is uniquely represented by a diagram of the form
  \eqref{equ:cartmor}.
\end{proposition}

\begin{proof}
    We have $A'\iso \poly{F'}(1)$ and $A\iso \poly{F}(1)$.
Define $\alpha : A' \rightarrow A$ to be the composite 
$$
\xymatrix{ A' \iso \poly{F'}(1) \ar[r]^{\phi_1} & \poly{F}(1) \iso A}\,.
$$
We need to construct $\beta:B' \to B$, and since it has to be
compatible with $\alpha$, $f'$ and $f$, it is enough to construct 
$B'_{a'} \to B_{\alpha(a')}$ for each $a'\in A'$.  Thereby we can 
reduce to the case where $A'=A=1$; in this case $\phi$ is invertible
since it is simultaneously vertical and cartesian.
But in this case the enriched 
Yoneda lemma above already ensures that the natural transformation is
induced by a unique map $B\to B'$, which we furthermore know is invertible. 
Its inverse is what we need for $B'_{a'} \to B_{\alpha(a')}$.
We have now constructed a diagram like \eqref{equ:cartmor}, and
it is routine to check that this diagram represents $\phi$.
\end{proof}

\begin{para}
  We give an example of a natural transformation that cannot be represented by
  diagrams.  On the category $\Set^{\mathbb{Z}_2}$ of involutive sets,
  the identity functor is represented by $1 \leftarrow 1 \to 1 \to 1$.  The
  twist natural transformation $\tau:\Id{}\Rightarrow\Id{}$, whose component on an
  object $X$ is the involution of $X$, is both cartesian and vertical.  It is
  clear that it cannot be represented by any diagram connecting $1\leftarrow 1
  \to 1\to 1$, since any connecting arrows would have to be identities and
  thereby induce the trivial natural transformation.  Observe that $\tau$ is
  not strong.
\end{para}

\begin{para} 
  We can now combine the diagrams representing vertical and cartesian strong
  natural transformations.  Given a diagram
  \begin{equation}
  \label{equ:morphism}
  \xycenter{
  G:&I \ar@{=}[dd]  & D  \ar[r]^{g} \ar[l]_{u} & C \ar[r]^{v} & J 
  \ar@{=}[dd] &  \\
   &  &  B'   \ar[r] \ar[d] \ar[u]  \drpullback & C   \ar@{=}[u] 
   \ar[d]   &&  \\
  F:&I  & B \ar[r]_{f} \ar[l]^{s} & A \ar[r]_{t}   &J &} 
  \end{equation}
  there is induced, by \ref{para:cart} and \ref{para:vert}, a strong natural
  transformation $\poly{\phi} : \poly{G} \Rightarrow \poly{F}$.  We refer to a
  diagram like~\eqref{equ:morphism} as a \myemph{morphism} from $G$ to $F$.  We
  arrive at the following result, a version of which appears as~\cite[Theorem
  3.4]{AbbottM:catc}, where it is stated for polynomial functors between slice
  categories over discrete objects.
\end{para} 

\begin{theorem}\label{thm:fibnat}
  Every strong natural transformation $\poly{G} \Rightarrow \poly{F}$ between
  polynomial functors is represented in an essentially unique way by a diagram like
  \eqref{equ:morphism}.
\end{theorem}

\begin{proof}
  By Proposition~\ref{thm:grothfib}, every strong natural transformation factors
  as a vertical strong transformation followed by a cartesian strong natural
  transformation in an essentially unique way.  The claim then follows from
  Proposition~\ref{thm:vertmor} and Proposition~\ref{thm:cartmor}.
\end{proof} 

\begin{corollary} 
  Every strong natural transformation between polynomial functors is a composite
  of units and counits of the basic adjunctions, their inverses when they
  exist, and coherence $2$-cells for pullback and
  its adjoints.
\end{corollary} 

\begin{proof} 
  The ingredients of the constructions in \ref{para:cart} and \ref{para:vert}
  are units, counits, pseudo-functoriality $2$-cells, as well as
  Beck-Chevalley isomorphisms, which in turn are
  constructed using units and counits (and inverses of their composites). 
\end{proof} 

\begin{para}\label{para:ff}
  Polynomials from $I$ to $J$ and their morphisms form a category denoted
  $\PolyNat(I,J)$.  Vertical composition of diagrams like~\eqref{equ:morphism}
  involves a simple pullback construction that via extension amounts precisely
  to refactoring cartesian-followed-by-vertical into
  vertical-followed-by-cartesian, cf.~the fibration property.  This can also be
  described as the unique way of defining vertical composition of diagrams to make
  the assignment
  given by extension functorial.  If we let $\PolyFun(\slice{I}, \slice{J})$ denote
  the category of polynomial functors from $\slice I$ to $\slice J$ and strong natural
  transformations, we can reformulate Theorem~\ref{thm:fibnat}
  as follows.
\end{para}

\begin{lemma}\label{thm:fij}
  For any $I,J$, the functor given by extension,
  $$
  \ext: \PolyNat(I,J) \rightarrow \PolyFun(\slice{I}, \slice{J})\,,
  $$
  is an equivalence of categories.
\end{lemma}
  
\begin{para}
  The involved categories are hom categories of appropriate bicategories of
  polynomials and polynomial functors, respectively, that we now describe,
  assembling the equivalences of the lemma into a biequivalence of bicategories
  (\ref{thm:bieq}).  We define the $2$-category of polynomial functors
  $\PolyFun$ as the sub-$2$-category of $\Cat$ having slices of $\catE$ as
  $0$-cells, polynomial functors as $1$-cells, and strong natural
  transformations as $2$-cells.

  We shall describe a bicategory $\PolyNat$ which has objects of $\catE$ as
  $0$-cells, polynomials as $1$-cells, and whose $2$-cells are the morphisms of
  polynomials, i.e.~diagrams like \eqref{equ:morphism}.  The vertical
  composition of $2$-cells has already been described, as has the horizontal
  composition of $1$-cells.  To define the horizontal composition of $2$-cells
  we simply transport back the $2$-cell structure from $\PolyFun$ along the
  local equivalences of Lemma~\ref{thm:fij}.

We begin by extending the family of functions mapping a pair of composable
polynomials $F$ and $G$ to their composite $G \circ F$, which we defined in
Paragraph~\ref{para:comp}, to a family of functors
\[
\PolyNat(J,K) \times \PolyNat(I,J) \to \PolyNat(I,K) \, .
\]
For this, let $\phi : F \Rightarrow F'$ be a morphism between polynomials from
$I$ to $J$, and let $\psi : G \Rightarrow G'$ be a morphism between polynomials
from $J$ to $K$.  We define the morphism $\psi \bcomp \phi : G \bcomp F
\Rightarrow G' \bcomp F'$ as the unique morphism of polynomials making the
following diagram commute
\[
\xymatrix{
\exppoly{G \bcomp F} \ar[rr]^{ \exppoly{\psi \bcomp \phi} }  
\ar[d]_{ \alpha_{G,F}  } & & \exppoly{ G' \bcomp F'  }  \ar[d]^{  \alpha_{G',F'}   } \\
\exppoly{G} \comp \exppoly{F}  \ar[rr]_{ \exppoly{\psi} \comp \exppoly{\phi}  } 
& & \exppoly{G'} \comp \exppoly{F'} \,. }
\]
Here $\alpha_{G,F}$ and $\alpha_{G',F'}$ are instances of
the isomorphism of Theorem~\ref{thm:subst},
and the diagram now expresses the naturality of $\alpha$.
We therefore get the following natural isomorphism of functors
\[
\xymatrix{
\PolyNat(J,K) \times \PolyNat(I,J) \ar[r] 
\ar@{}[dr] | (.59){\iso} 
\ar[d]_{ P_{J,K} \times P_{I,J}} & 
\PolyNat(I,K) 
\ar[d]^{ P_{I,K}} \\
\PolyFun(\slice{J},\slice{K}) \times 
\PolyFun(\slice{I},\slice{J}) \ar[r]  & 
\PolyFun(\slice{I},\slice{K}) }
\]
where the top horizontal functor is substitution of polynomials and the
bottom horizontal map is composition of functors in $\PolyFun$.
The identity maps in $\PolyNat$ are represented by the polynomials 
$\id{I} : I \rightarrow I$, and  we have natural isomorphisms
\[
\xymatrix{
                    & \PolyNat(I,I) \ar[dd]^{P_{I,I}} \\
\defcat{1} \ar[ur]^(.45){\id{I}}  \ar[dr]_(.45){1_{\slice{I}}}  \ar@{}[r] | (.55){\iso} &   \\
                    & \PolyFun(\slice{I},\slice{I})  \,.}
\]
We define the associativity and unit isomorphisms. For associativity, 
given polynomials  $F : I \rightarrow J$, $G : J \rightarrow K$, and $H : K 
\rightarrow L$,
define 
\[
\theta_{H,G,F} : (H \bcomp G) \bcomp F \Rightarrow H \bcomp (G \bcomp F)
\]
to be the unique morphism of polynomials making the following diagram commute
\begin{equation}
\label{equ:psdfun1}
\xycenter{
\exppoly{(H \bcomp G) \bcomp F} 
\ar[r]^{\exppoly{\theta_{H,G,F}}} 
\ar[d]_{\alpha_{H \bcomp G, F}} & 
\exppoly{H \bcomp (G \bcomp F) } \ar[d]^{\alpha_{H, G \bcomp F}} \\
\exppoly{H \bcomp G} \comp \exppoly{F} \ar[d]_{\alpha_{H,G} \comp 
\exppoly{F}} 
&  \exppoly{H} \comp \exppoly{G \bcomp F} \ar[d]^{\exppoly{H} \comp \alpha_{G,F}} \\
\big( \exppoly{H} \comp \exppoly{G} \big) \comp \exppoly{F} \ar@{=}[r] & 
\exppoly{H} \comp \big( \exppoly{G} \comp \exppoly{F} \big) \,. } 
\end{equation} 
For the unit isomorphisms, given a polynomial  $F : I \rightarrow J$, define
\[
\lambda_F : \id{J} \bcomp F \Rightarrow F \, , \qquad 
\rho_F: F \bcomp \id{I} \Rightarrow F 
\]
to be the unique morphism of polynomials such that
\begin{equation}
\label{equ:psdfun2}
\xycenter{
\exppoly{\id{J} \bcomp F } \ar[rr]^{\exppoly{\lambda_F}} \ar[d]_{\alpha_{\id{J}, F} }  & & 
\exppoly{F} \ar@{=}[d] \\
\exppoly{\id{J}} \comp \exppoly{F} \ar[rr]_{\alpha_J \comp \exppoly{F} }
& & 1_{\slice{J}} \comp \exppoly{F} }
\end{equation} 
and
\begin{equation}
\label{equ:psdfun3}
\xycenter{
\exppoly{F \bcomp \id{I}}  \ar[rr]^{\exppoly{\rho_F}} 
\ar[d]_{\phi_{F, \id{I} } }& & \exppoly{F} 
\ar@{=}[d] \\
\exppoly{F} \comp \exppoly{ \id{I} } \ar[rr]_{\exppoly{F} \comp \alpha_I }
& & \exppoly{F} \comp 1_{\slice{I}     }}
\end{equation} 
commute. All the data of the bicategory $\PolyNat$ have now been given.
The naturality and coherence axioms for a bicategory
can be verified by standard diagram-chasing arguments,
which exploit the uniqueness properties used to 
define the components of $\theta$, $\lambda$, and $\rho$.
The interchange law of $\PolyFun$ is used at several points. 
Let us remark that the definition of the bicategory $\PolyNat$ is essentially 
determined  by the requirement that we obtain a pseudo-functor 
\[
\ext : \PolyNat \rightarrow \PolyFun \, .
\]
Indeed, the diagrams in~\eqref{equ:psdfun1},~\eqref{equ:psdfun2},~\eqref{equ:psdfun3} 
express exactly the coherence conditions for a pseudo-functor~\cite{BenabouJ:intb}.
It is clear by construction that this pseudo-functor is bijective 
on objects, and it is locally an equivalence of categories by Lemma~\ref{thm:fij}.
Hence we have established the following.
\end{para}

\begin{theorem}\label{thm:bieq}
  The extension pseudo-functor
  \[
  \ext: \PolyNat \to \PolyFun
  \]
  is a biequivalence. \qed
\end{theorem}

\begin{para}\label{para:containers}
  The notions of polynomial and polynomial functors are almost exactly the same
  as what is called \myemph{container} and \myemph{container functor} by Abbott,
  Altenkirch and Ghani~\cite{AbbottM:phd, AbbottM:catc,
  Abbott-Altenkirch-Ghani:nested, Abbott-Altenkirch-Ghani:strictly-positive}.
  One minor difference is that they only consider slices over discrete objects,
  i.e.~of the form $\slice n \simeq \catE^n$, where $n$ denotes the sum of $n$
  copies of the terminal object, and for this they also need to assume finite 
  sums.  In our setting there is no reason for that
  restriction, and in fact Altenkirch and Morris~\cite{Morris-Altenkirch:lics09}
  have been able to lift the restriction also from the container theory,
  introducing the notion of \myemph{indexed container}.  Another difference, also
  quite minor, is that while we prefer to work with strength, the container
  people work with fibred categories, fibred functors and fibred natural
  transformations.  This involves replacing all slice categories $\slice I$ by
  the fibration over $\catE$ whose $K$-fibre is $\slice{(K\times I)}$, and work
  with those instead.  The two viewpoints are in fact equivalent, thanks to a
  result of Par\'e, who showed (cf.~\cite{JohnstoneP:carmt}) that if a strong
  functor preserves pullbacks then it is canonically indexed, i.e.~fibred.  (It
  is easy to see that a fibred functor has a strength.)  We have chosen the
  viewpoint of tensorial strength for its simplicity.  Modulo the above minor
  differences (and modulo Par\'e's theorem), Lemma~\ref{thm:carttopoly},
  Theorem~\ref{thm:fibnat}, and Theorem~\ref{thm:bieq} were also proved in
  Abbott's thesis~\cite{AbbottM:phd}.
%
\end{para}
  
\section{The double category of polynomial functors} 
\label{sec:double}

\begin{para}
  It is important to be able to compare polynomial functors with different
  endpoints, and to base change polynomial functors along maps in $\catE$.  This
  need can been seen already for linear functors~\ref{ex:linear}:
  a small category is a monad in
  the bicategory of spans \cite{BenabouJ:intb}, but in order to get functors
  between categories with different object sets, one needs maps between spans
  with different endpoints \cite{LackS:fortm2}.  The most convenient framework 
  for this is that of double categories, as it allows for diagrammatic
  representation.  The base change structure is concisely captured in Shulman's
  notion of framed bicategory~\cite{ShulmanM:frabmf}: our double categories of
  polynomial functors will in fact be framed bicategories.
\end{para}

\begin{para}
  Recall that a double category $\mathbb{D}$ consists of a category of objects
  $\mathbb{D}_0$, a category of morphisms $\mathbb{D}_1$, together with
  structure functors
    \[
      \xymatrix @+5pt {
    \mathbb{D}_0 \ar[r] & \ar@<.9ex>[l]^{\partial_0} \ar@<-.9ex>[l]_{\partial_1} \mathbb{D}_1 &  
    \ar[l]_-{\text{\tiny{comp.}}} \mathbb{D}_1 \times_{\mathbb{D}_0} \mathbb{D}_1 }
      \]
  subject to the usual category axioms.  The objects of $\mathbb{D}_0$ are
  called \myemph{objects} of $\mathbb{D}$, the morphisms of $\mathbb{D}_0$ are
  called \myemph{vertical arrows}, the objects of $\mathbb{D}_1$ are called
  \myemph{horizontal arrows}, and the morphisms of $\mathbb{D}_1$ are called
  \myemph{squares}.  As is custom \cite{GrandisM:adjdc}, we allow the
  possibility for the horizontal composition to be associative and unital only
  up to specified coherent isomorphisms.  Precisely, a double category is a
  pseudo-category~\cite{MartinsFerreiraN:psec} in the $2$-category $\Cat$; see
  also \cite[\S 5.2]{LeinsterT:higohc}.
\end{para}

\begin{para} 
  A \myemph{framed bicategory} \cite{ShulmanM:frabmf} is a double category for which 
  the functor
 \[
  (\partial_0,\partial_1): \mathbb{D}_1 \longrightarrow \mathbb{D}_0\times\mathbb{D}_0
  \]
  is a bifibration.  
  (In fact, if it is a fibration then it is automatically an opfibration, 
  and vice versa.) 
  The upshot of this condition is that horizontal arrows can be base changed and
  cobase changed along arrows in $\mathbb{D}_0\times\mathbb{D}_0$ (i.e.~pairs
  of vertical arrows).
\end{para}

\begin{para} 
  We need to fix some terminology.  The characteristic property of a fibration
  is that every arrow in the base category admits a cartesian lift, and that
  every arrow in the total space factors (essentially uniquely) as a vertical
  arrow followed by a cartesian one.  In the present situation, the term
  `cartesian' is already in use to designate cartesian natural transformations
  (which fibrationally speaking are vertical rather than cartesian), and the
  word `vertical' already has a double-categorical meaning.  For these reasons,
  instead of talking about `cartesian arrow' for a fibration we shall say
  \myemph{transporter} arrow; this terminology goes back to
  Grothendieck~\cite{GrothendieckA:revegf}.  Correspondingly, we shall say
  \myemph{cotransporter} instead of opcartesian.  We shall simply refrain from
  using `vertical' in the fibration sense.  The arrows mapping to identity
  arrows by the fibration will be precisely the natural transformations of
  polynomial functors.
\end{para} 

\begin{para}
  We want to extend the bicategories $\PolyNat$ and $\PolyFun$ to double
  categories.  The objects of the double category $\DblPolyFun$ are the slices
  of $\catE$, and the horizontal arrows are the polynomial functors.  The
  vertical arrows are the dependent sum functors (i.e.~functors of the form
  $\ladj{u}$ for some $u$), and the squares in $\DblPolyNat$ are of the form
  \begin{equation}\label{sw}
  \xycenter 
  {
  \slice{I'} \ar[d]_{\ladj{u}} \ar[r]^{P'} \ar@{}[dr]|{\ \Downarrow \, \phi} 
  & 
  \slice{J'} \ar[d]^{\ladj{v}} \\
  \slice{I} \ar[r]_{P} & \slice{J}
  }
  \end{equation}
  where $P'$ and $P$  are polynomial functors and $\phi$  is a strong natural 
  transformation. 
\end{para}

\begin{proposition}\label{thm:frabic}
  The double category $\DblPolyFun$ is a framed bicategory.
\end{proposition}

\begin{proof}
  The claim is that the functor sending a polynomial functor $P : \slice{I}
  \rightarrow \slice{J}$ to~$(I, J)$ is a bifibration.  For each pair of arrows
  $u : I' \rightarrow I$, $v : J' \rightarrow J$ in $\catE$ we have the
  following basic squares (companion pairs and conjoint
  pairs~\cite{GrandisM:adjdc})
\begin{equation*}
  \xymatrixrowsep{50pt}
  \xymatrixcolsep{50pt}
  \xymatrix @!=0pt {
  \slice{I'}  \ar[r]^{\ladj{u}} \ar[d]_{\ladj{u}} 
  & \slice{I}  \ar@{=}[d] \\
  \slice{I}  \ar@{=}[r]  & \slice{I}   }
  \quad
  \xymatrix @!=0pt {
  \slice{I'}  \ar@{=}[r] \ar[d]_{\ladj{u}} \ar@{}[dr]|{\  \Downarrow \,  \eta}
  &  \slice{I'}  \ar@{=}[d] \\
  \slice{I} \ar[r]_{\pbk{u}}  & \slice{I'}    }
\quad
  \xymatrix @!=0pt {
  \slice{J'}  \ar@{=}[r] \ar@{=}[d] 
  & \slice{J'}  \ar[d]^{\ladj{v}} \\
  \slice{J'} \ar[r]_{\ladj{v}}  & \slice{J} 
  }
  \quad
  \xymatrix @!=0pt {
  \slice{J} \ar[r]^{\pbk{v}} \ar@{=}[d] \ar@{}[dr]|{\ \Downarrow \, \varepsilon} 
  &
  \slice{J'}  \ar[d]^{\ladj{v}} \\
  \slice{J}  \ar@{=}[r]  & \slice{J} 
  }
\end{equation*}
It is now direct to check that the pasted square
\[
  \xymatrixrowsep{48pt}
  \xymatrixcolsep{56pt}
  \xymatrix @!=0pt {
  \slice{I'} \ar[r]^{\ladj{u}} \ar[d]_{\ladj{u}} 
   & 
  \slice{I} \ar@{=}[d] 
  \ar[r]^P & \slice{J}\ar@{=}[d] 
  \ar[r]^{\pbk{v}}  \ar@{}[dr]|{\ \Downarrow \, \varepsilon} 
  & \slice{J'} \ar[d]^{\ladj{v}}\\
  \slice{I} \ar@{=}[r]  & \slice{I} \ar[r]_P & \slice{J} \ar@{=}[r] & \slice{J}
  }
\]
is a transporter lift  of $(u, v)$  to $P$. We call $\pbk{v} \circ P \circ \ladj{u}$ the \myemph{base change}
   of $P$ along $(u,v)$, and denote it   $\basechange{(u,v)} (P)$.

   Dually,  it is direct to check that the pasted square
\[
  \xymatrixrowsep{48pt}
  \xymatrixcolsep{56pt}
  \xymatrix @!=0pt {
  \slice{I'} \ar@{=}[r] \ar[d]_{\ladj{u}} \ar@{}[dr]|{\ \Downarrow \, \eta} 
  & 
  \slice{I'} \ar@{=}[d] 
  \ar[r]^{P'} & \slice{J'}
  \ar@{=}[d] 
  \ar@{=}[r] 
  & \slice{J'} \ar[d]^{\ladj{v}}\\
  \slice{I} \ar[r]_{\pbk{u}}  & \slice{I'} \ar[r]_{P'} &
  \slice{J'} \ar[r]_{\ladj{v}} & \slice{J}
  }
\]
is a cotransporter lift  of $(u, v)$ to $P'$. We call $\ladj{v} \circ P' \circ \pbk{u}$
the \myemph{cobase change}
of $P'$ along $(u,v)$, and denote it $\cobasechange{(u,v)}(P')$.
\end{proof}
The above procedure of getting a framed bicategory out of a bicategory is a
general construction: one starts with a bicategory $\catC$ with a subcategory
$\catL$ consisting of left adjoints and comprising all the objects of $\catC$,
and obtains a framed bicategory by taking as vertical arrows the arrows in
$\catL$.  The details can be found in \cite[Appendix]{ShulmanM:frabmf}.

\begin{para}
  Via the biequivalence $\PolyNat\simeq \PolyFun$ 
  between the bicategory of polynomials and the
  $2$-category of polynomial functors, Proposition~\ref{thm:frabic} gives us
  also a framed bicategory of polynomials $\DblPolyNat$, featuring nice diagrammatic
  representations which we now spell out, extending the results of
  Section~\ref{sec:morphisms}.  The following result is
  the double-category version of Theorem~\ref{thm:fibnat}.
\end{para}

\begin{theorem}
  The squares \textrm{(\ref{sw})} of $\DblPolyFun$ are represented by diagrams of the
  form
\begin{equation}
\label{equ:sqdiag} 
  \xycenter{ 
P':&  I'\ar[dd]_u & \ar[l]  B' \ar[r]& A' \ar@{=}[d] \ar[r]&J'\ar[dd]^v&\\
&  &\cdot \drpullback\ar[u] \ar[r] \ar[d] &\cdot \ar[d]&& \\
P:&  I & \ar[l] B \ar[r]& A \ar[r]&J \, .&\\
  }
\end{equation} 
This representation is unique up to choice of pullback in the middle. 
It follows that extension constitutes a framed biequivalence 
\[
\DblPolyNat \isopil \DblPolyFun \,.
\]
\end{theorem}

\begin{proof}
  By Theorem~\ref{thm:fibnat}, diagrams like~(\ref{equ:sqdiag}) (up to 
  choice of pullback) are in
  bijective correspondence with strong natural transformations $\ladj{v} \circ P' \circ
  \pbk{u} \Rightarrow P$, which by adjointness correspond to strong natural
  transformations $\ladj{v} \circ P' \Rightarrow P\circ \ladj{u}$, 
  i.e.~squares~(\ref{sw})
  in $\DblPolyFun$.
\end{proof}

\begin{para}
  The vertical composition of two diagrams
  \[
  \xycenter{ 
\cdot\ar[dd] & \ar[l]  \cdot \ar[r]& \cdot \ar@{=}[d] \ar[r]&\cdot\ar[dd]\\
&\cdot \drpullback\ar[u] \ar[r] \ar[d] &\cdot \ar[d]& \\
\cdot\ar[dd] & \ar[l]  \cdot \ar[r]& \cdot \ar@{=}[d] \ar[r]&\cdot\ar[dd]\\
&\cdot \drpullback\ar[u] \ar[r] \ar[d] &\cdot \ar[d]&\\
\cdot & \ar[l] \cdot \ar[r]& \cdot \ar[r]&\cdot\\
  }
\]
  is performed by replacing the two middle squares
   \[
  \xycenter{ 
\cdot \drpullback \ar[r] \ar[d] &\cdot \ar[d]\\
\cdot \ar[r] & \cdot \ar@{=}[d]\\
\cdot \ar[u] \ar[r] &\cdot
  }
\]
by a configuration 
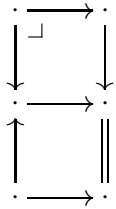
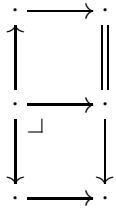
   \[
  \xycenter{ 
\cdot \ar[r] & \cdot \ar@{=}[d]\\
\cdot \drpullback \ar[r] \ar[u]\ar[d] &\cdot \ar[d]\\
\cdot \ar[r] &\cdot
}
\]
and then composing vertically.  The replacement is a simple pullback 
construction, and
checking that the composed diagram has the same extension as the
vertical pasting of the extensions is a straightforward calculation.
\end{para}

\begin{para}
  At the level of polynomials, the bifibration $\DblPolyNat \to 
  \catE\times\catE$ is now the `endpoints' functor, associating
  to a polynomial
  $I \leftarrow B \to A \to J$ the pair $(I,J)$.
  With notation as in the proof of Proposition~\ref{thm:frabic},
  we know the cobase change of $P'$ along $(u,v)$ is just $\ladj{v} \circ P'
  \circ \pbk{u}$, and it is easy to see that
\[
\xymatrixrowsep{36pt}
 \xymatrixcolsep{45pt}
 \xymatrix @!0 {
P':  & I'\ar[d]_u & \ar[l] B'  \ar[r]   \ar@{=}[d]& A' \ar@{=}[d] 
\ar[r]&J'\ar[d]^v &\\
\cobasechange{(u,v)}(P'):& I & \ar[l] B'  \ar[r]& A' \ar[r]&J &
 }
 \]
is a cotransporter lift of $(u,v)$ to $P'$.

The transporter lift of $(u,v)$ to $P$, which is the same thing as the base
change of $P$ along $(u,v)$, is slightly more complicated to construct.  It can
be given by first base changing along $(u,\id{})$ and then along
$(\id{},v)$:
\begin{equation}\label{eq:basechange}
\xymatrixrowsep{32pt}
\xymatrixcolsep{45pt}
\xymatrix @!0 {
\basechange{(u,v)}(P): &  I'\ar@{=}[d] & \ar[l] \cdot  \drpullback \ar[r] 
\ar[d]& \cdot \drpullback\ar[d] \ar[r]&J'\ar[d]^v &\\
\basechange{(u,\id{})}(P):&\cdot \ar@{=}[d] & \cdot \drpullback \ar[l] \ar[d] \ar[r] & \cdot 
\ar[dd] \ar[r] & \cdot \ar@{=}[dd] &\\
&\cdot \ar[d]_u & \cdot \dlpullback \ar[l] \ar[d] && &\\
P : &  I & \ar[l] B  \ar[r]& A \ar[r]&J 
}
\end{equation}
\end{para}

\begin{para}\label{sourcelift}
The intermediate polynomial $\basechange{(u,\Id{})}(P)$ is called the 
\myemph{source lift} of $P$ along $u$, and we shall need it later on.
Since $\partial_0$ (as well as $\partial_1$) is itself a bifibration,
for which the source lift is the transporter lift, it
enjoys the following universal property: every 
square 
\[
  \xymatrixrowsep{32pt}
  \xymatrixcolsep{42pt}
  \xymatrix @!0 {
P':&  I'\ar[dd]_u & \ar[l]  \cdot \ar[r]& \cdot \ar@{=}[d] \ar[r]&J'\ar[dd]^v&\\
&  &\cdot \drpullback\ar[u] \ar[r] \ar[d] &\cdot \ar[d]&& \\
P:&  I & \ar[l] \cdot \ar[r]& \cdot\ar[r]&J &\\
  }
\]
factors uniquely through the source lift, like 
\[
  \xymatrixrowsep{32pt}
  \xymatrixcolsep{42pt}
  \xymatrix @!0 {
 P':&   I'\ar@{=}[dd] & \ar[l] \cdot \ar[r]&\cdot\ar@{=}[d] \ar[r]&J'\ar[dd]^v &\\
  &  &\cdot \drpullback\ar[u] \ar[r] \ar[d] &\cdot \ar[d]& &\\
\basechange{(u,\id{})}(P):&I' \ar[d]_u & \cdot \drpullback \ar[l] \ar[d] \ar[r] & \cdot 
\ar[d] \ar[r] & \cdot \ar@{=}[d] &\\
P : &  I & \ar[l] \cdot  \ar[r]& \cdot \ar[r]&J &
}
\]
where the bottom part is as in (\ref{eq:basechange}).
\end{para}

\begin{para}
  All the constructions and arguments of this section apply equally well inside the
  cartesian fragment: starting with the $2$-category $\PolyFunCart$
  of polynomial functors and
  their cartesian strong natural transformations, a double category
  $\DblPolyFunCart$ results, which is a framed bicategory.  The only point to note
  is that all the constructions are compatible with the cartesian condition,
  since they all depend on the $\ladj{}\adjoint\pbk{}$ adjunction, which is
  cartesian.  Note also that the transporter and cotransporter lifts belong to
  the cartesian fragment.  The following two results follow readily.
\end{para}

\begin{proposition}
  The double category $\DblPolyFunCart$ whose objects are the slices of
  $\catE$, whose horizontal arrows are the polynomial functors, whose
  vertical arrows are the dependent sum functors, and whose squares are 
  cartesian strong natural transformations
\[
  \xycenter{
  \slice{I'} \ar[d]_{\ladj{u}} \ar[r]^{P'} \ar@{}[dr]|{\ \Downarrow \, \phi} 
 & 
  \slice{J'} \ar[d]^{\ladj{v}} \\
  \slice{I} \ar[r]_{P} & \slice{J}
  }
\]
  is a framed bicategory.
\end{proposition}

\begin{proposition}
  The squares of $\DblPolyFunCart$ are represented uniquely by diagrams
\begin{equation}
\label{cartesian2cell}  
\xycenter{ 
  I'\ar[d] & \ar[l] B'  \drpullback\ar[r] \ar[d]& A' \ar[d] \ar[r]&J'\ar[d] \\
  I & \ar[l] B  \ar[r] & A \ar[r]&J \,, \\
}
\end{equation}
hence extension constitutes a framed biequivalence
\[
\DblPolyCart \isopil \DblPolyFunCart \,.
\]
\end{proposition}

\begin{para}
For the remainder of this paper, we shall only deal with the cartesian fragment,
which is also what is needed in \cite{KockJ:polft} and
\cite{KockJ:polfo}.  In those two papers, a central construction is to label
trees by a polynomial endofunctor $P$.  Trees are themselves seen as polynomial
endofunctors (cf.~Example~\ref{ex:trees}), and the labelling amounts precisely 
to a cartesian $2$-cell in the double category of polynomial functors.  The
importance of the cartesian condition (a bijection of certain fibres) is to
ensure that a node in a tree is labelled by an operation of the same arity.
\end{para}

\begin{para}\label{para:shaf}
  We finish this section with a digression on the relationship between
  polynomial functors and the shapely functors and shapely types of
  Jay and Cockett~\cite{JayCB:shatp,JayCB:sems}, since the double-category
  setting provides some conceptual simplification of the latter notion.
  We now assume $\catE$ has sums.
  
  A \myemph{shapely functor} \cite{JayCB:shatp} is a pullback-preserving functor
  $F : \catE^m \rightarrow \catE^n$ equipped with a strength.  Since, for a
  natural number $n$, the discrete power $\catE^n$ is equivalent to the slice
  $\catE/n$, where $n$ now denotes the $n$-fold sum of $1$ in $\catE$, it makes
  sense to compare shapely functors and polynomial functors.  Since a polynomial
  functor preserves pullbacks and has a canonical strength, it is canonically a
  shapely functor.  It is not true that every shapely functor is polynomial.
  For a counter example, let $K$ be a set with a non-principal filter $\catD$,
  and consider the filter-power functor
  \begin{eqnarray*}
    F:\Set & \longrightarrow & \Set\\
    X & \longmapsto & \colim_{D\in \catD} X^D \,,
  \end{eqnarray*}
  which preserves finite limits since it is a filtered colimit of
  representables.  Since every endofunctor on $\Set$ has a canonical strength,
  $F$ is a shapely functor.  However, $F$ does not preserve all cofiltered
  limits, and hence, by \ref{para:six}\,(\ref{item:connected}) cannot be
  polynomial.  For example, $\emptyset = \lim_{D\in\catD} D$ itself is not
  preserved.  This example is apparently at odds with Theorem~8.3 of
  \cite{AbbottM:catc}.
\end{para} 

\begin{para}\label{para:shat}
  Let $L: \catE \to \catE$ denote the \myemph{list endofunctor},
  $L(X)=\sum_{n\in\Nat}X^n$, which is the same as what we called the
  free-monoid monad in Example~\ref{ex:freemonoid}.  A \myemph{shapely 
  type}~\cite{JayCB:shatp} in
  one variable is a shapely functor equipped with a cartesian strong natural
  transformation to $L$.  A morphism of shapely types is a natural
  transformation commuting with the structure map to~$L$.  The idea is that the
  shapely functor represents the template or the shape into which some data can
  be inserted, while the list holds the actual data; the cartesian natural
  transformation encodes how the data is to be inserted into the template.  As
  emphasized in~\cite{MoggiE:monsft}, the cartesian strong natural
  transformation is part of the structure of a shapely type.
  Since any functor with a cartesian natural transformation to $L$ is polynomial
  by Lemma~\ref{thm:carttopoly}, it is clear that one-variable shapely types are essentially
  the same thing as one-variable polynomial endofunctors with a cartesian natural
  transformation to $L$, and that there is an equivalence of categories between
  the category of shapely types and the category $\PolyCart(1,1)/L$.
  
  According to Jay and Cockett~\cite{JayCB:shatp}, a shapely type in $m$ input 
  variables and $n$ output variables is a shapely functor $\catE^m \to \catE^n$
  equipped with a cartesian strong natural transformation to the functor $L_{m, n} :
  \catE^m \rightarrow \catE^n$ defined by
\[
\textstyle
L_{m,n}(X_i \mid  i \in m) = \big( L(\sum_{i \in m} X_i) \mid j \in n \big)\,,
\]
and they motivate this definition by considerations on how to insert data into 
templates.
With the double-category formalism, we can give a conceptual explanation
of the formula:
writing $u_m : m\to 1$ and $u_n : n \to 1$  for the maps to the terminal object, the 
functor $L_{m, n} : \catE^m \rightarrow \catE^n$ is nothing but the composite
\[
\pbk{u_n} \circ L \circ \ladj{u_m} = \basechange{(u_{m},u_{n})} L \,,
\]
the base change of $L$ along $(u_m,u_n)$.
Hence we can say uniformly that a
shapely type is an object in $\PolyCart/L$ with endpoints
finite discrete objects.  
\end{para}

\section{Polynomial monads} 
\label{sec:polmnd}

\begin{para} \label{para:pmonad}
  Let $I \in \catE$.  A \myemph{polynomial monad} on $\slice I$ is a monad $(T,
  \eta, \mu)$ for which $T$ is a polynomial functor and $\eta$ and $\mu$ are
  cartesian strong natural transformations.  From the point of view of the
  formal theory of monads~\cite{StreetR:fortm}, a polynomial monad is a monad in
  the $2$-category~$\PolyFunCart$.  A basic example of a polynomial monad is
  the free-monoid monad of Example~\ref{ex:freemonoid}.
\end{para}

\begin{para}
  We are interested in the construction of the free monad on a polynomial
  endofunctor, and start by recalling from~\cite{KellyG:unittc,BarrM:toptt} some
  general facts about free monads.  Let $\catC$ be a category and $P : \catC
  \rightarrow \catC$ an endofunctor.  The \myemph{free monad} on $P$ is a monad
  $(T,\eta,\mu)$ on $\catC$ together with a natural transformation $\alpha : P
  \Rightarrow T$ enjoying the following universal property: for any monad $(T',
  \eta', \mu')$ on $\catC$ and any natural transformation $\phi : P \Rightarrow
  T'$ there exists a unique monad morphism $\phi^\sharp : T \Rightarrow T'$
  making the following diagram commute:
  \[
  \xymatrix{
  P \ar[r]^{\alpha} \ar@/_1pc/[dr]_{\phi} & T \ar[d]^{\phi^\sharp} \\
		  & T'}
  \]
  The following construction of the free monad on $P$ is standard.
  Let $\Palg$ denote the category of
  $P$-algebras and $P$-algebra morphisms.  We denote $P$-algebras as pairs
  $(X,\mysup_X)$ where $X$ is the underlying object, and $\mysup_X: PX \to X$ 
  is the structure map, sometimes suppressed from the notation for brevity.
  If the forgetful functor $U :
  \Palg \rightarrow \catC$ has a left adjoint, then the monad $(T, \eta, \mu)$
  resulting from the adjunction is the free monad on $P$.
  If $\catC$ has binary sums, a necessary and sufficient condition for the
  existence of the left adjoint to $U$
is that, for every $X \in \catC$, the endofunctor $X + P(-) : \catC 
\rightarrow \catC$ has an initial algebra.  Indeed, in that case we can construct
the free monad as follows. For $X \in \catC$, we define
$TX$ as the initial algebra for $X + P(-) : \catC \rightarrow \catC$,
and $\eta_X : X \rightarrow TX$ as the composite
\[
\xymatrix{
X \ar[r]^(.3){\iota_1} & X + P(TX) \ar[r]^(.65){t_X} & TX }
\]
where $\iota_1$ is the first sum inclusion and $t_X$ is the structure
map of $TX$ as an $(X+P)$-algebra.  Finally, since $T^2 X$ is the initial algebra for the
functor $TX + P(-)$, we can define $\mu_X : T^2 X \rightarrow TX$ 
as the unique map making the following diagram commute:
\[
\xymatrix{
TX + P(T^2X) \ar[rr]^{TX + P(\mu_X)} \ar[dd]_{t_{TX}} & &   TX + P(TX) \ar[d] \\
                            &  & TX + X + P(TX) \ar[d]^{(1_{TX}, t_X)} \\
T^2 X \ar[rr]_{\mu_X}                &  & TX.} 
\]
Functoriality, naturality, and the monad axioms follow readily from these
definitions.  Note that the $X$-component of the natural transformation $\alpha : P
\Rightarrow T$ is given as the composite
\begin{equation}
  \label{alpha}
\xymatrix{
PX \ar[r]^(.4){P(\eta_X)} & P(TX) \ar[r]^(.42){\iota_2} & X + P(TX) 
\ar[r]^(.6){t_X} & TX \,. }
\end{equation}
\end{para}

\begin{para}\label{para:W}
  Let us now return to the locally cartesian closed category $\catE$, now 
  assumed to be extensive and in particular have finite sums.
  Recall from~\cite{MoerdijkI:weltc} that $\catE$ is said to have W-types if
  every polynomial functor in a single variable on $\catE$ has an initial
  algebra.  This terminology is motivated by the fact that initial algebras for
  polynomial functors in a single variable are category-theoretic counterparts
  of Martin-L\"of's types of wellfounded trees~\cite{NordstromB:marltt}.  Every
  elementary topos with a natural numbers object has
  W-types~\cite{MoerdijkI:weltc}.  If $\catE$ has W-types, then every polynomial
  endofunctor, not just those in a single variable, has an initial
  algebra~\cite[Theorem 14]{GambinoN:weltdp}.  Initial algebras for general
  polynomial functors are category-theoretic counterparts of Petersson and
  Synek's general tree types~\cite{PeterssonK:setcis}; see
  also~\cite[Chapter~16]{NordstromB:promlt}.

  Henceforth, we assume that $\catE$ has W-types.  For any polynomial
  endofunctor $P : \slice{I} \rightarrow \slice{I}$ and any $X\in \slice I$, the
  functors $X + P(-) : \slice{I} \rightarrow \slice{I}$ are again polynomial,
  hence have initial algebras.  Therefore every polynomial endofunctor admits a
  free monad.
\end{para}

\begin{para}
  Theorem~\ref{thm:freemonad} below asserts that the free monad on a polynomial
  functor is polynomial.  The proof exploits the possibility of recursively
  defining maps out of initial algebras for polynomial functors, and we need 
  first to set up some notation to handle this.
  Let $P : \slice{I} 
\rightarrow \slice{I}$ be the polynomial functor represented by the 
diagram
\[
\xycenter{
I  & \ar[l]_{s} B \ar[r]^{f}  & A \ar[r]^{t} & I \,. }
\]
We regard such a diagram as a generalised many-sorted signature.  This point of
view is most easily illustrated by considering the case of~$\catE = \Set$.  The
object $I$ provides the set of sorts of the signature.  The set of terms of the
signature is defined inductively by saying that we have a term
$\mathrm{sup}_a(x)$ of sort $t(a)$ whenever $a \in A$ and $x = (x_b \mid b \in
B_a)$ is a family of terms such that $x_b$ has sort $s(b)$ for all~$b \in B_a$.
Such a term may be represented graphically as a one-node tree
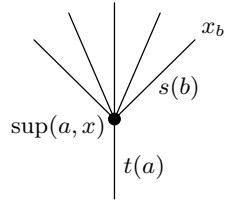
\begin{center}\begin{texdraw}
  \move (0 0) \lvec (0 30) \fcir f:0 r:2.5
  \lvec (-30 60)
  \move (0 30) \lvec (-17 70)
  \move (0 30) \lvec (0 74)
  \move (0 30) \lvec (17 70)
  \move (0 30) \lvec (30 60)
  
  \htext (11 12) {\footnotesize $t(a)$}
  \htext (24 42) {\footnotesize $s(b)$}
  \htext (37 64) {\footnotesize $x_b$}
  \htext (-21 27) {\footnotesize $\mysup(a,x)$}
\end{texdraw}\end{center}
The incoming edges are indexed by the elements of $B_a$ and further labelled by
elements of $I$, with the edge indexed by $b \in B_a$ labelled by $s(b) \in I$.
The outgoing edge is labelled by $t(a) \in I$.  We label the node $\mysup(a,x)$ if
the family $x = (x_b \mid b \in B_a)$ labels its incoming edges.

  Let $W$ be the initial algebra for $P$, with structure map $\mysup_W : PW
  \rightarrow W$.  Initiality of the algebra means that for any other algebra
  $(X,\mysup_X)$, there exists a unique algebra map $\theta : W \rightarrow X$,
  thus making the following diagram commute
\[
\xymatrix{
PW \ar[r]^{P(\theta)} \ar[d]_{\mysup_W} & PX \ar[d]^{\mysup_X} \\
W \ar[r]_{\theta} & X \,.}
\]
In the internal language of $\catE$, we can represent the structure
map of $W$ as the $I$-indexed family
\[
\mysup_{W_i} : \sum_{a \in A_i} \prod_{b \in B_a} W_{sb} \rightarrow W_i \, .
\]
The initiality of $W$ can be expressed by
saying that there exists a unique family of maps $\theta_i :
W_i \rightarrow X_i$  satisfying the recursive equation
\[
\theta_i(\mysup_{W_i}(a,h)) = \mysup_{X_i}(a, (\lambda b \in B_a) \; \theta_{sb}( hb))
\,,
\]
where we employ lambda notation $(\lambda b \in B_a) \; \theta_{sb}( hb)$ to 
indicate the function $B_a \to X$ sending $b$ to $\theta_{sb}(hb)$.
\end{para}

\begin{theorem}\label{thm:freemonad}
  The free monad on a polynomial endofunctor is a polynomial monad.
\end{theorem}

\begin{proof}
  Let $P : \slice{I} \rightarrow \slice{I}$ be the polynomial endofunctor
  represented by
  \[
  \xycenter{ I  & \ar[l]_{s} B \ar[r]^{f}  & A \ar[r]^{t} & I \,, }
  \]
  and let $(T, \eta, \mu)$ be the free monad on $P$.  We need to show that $T :
  \slice{I} \rightarrow \slice{I}$ is a polynomial functor, and that $\eta :
  \Id{} \Rightarrow T$ and $\mu : T^2 \Rightarrow T$ are cartesian strong
  natural transformations.  We shall show that $T$ is naturally isomorphic to
  the polynomial functor represented by the diagram
  \begin{equation}\label{equ:DC}
  \xymatrix{
  I  & \ar[l]_{u} D \ar[r]^g  & C \ar[r]^v & I }
  \end{equation}
  whose constituents we now proceed to construct.  Intuitively, $C$ is the set
  of wellfounded trees with branching profile given by the polynomial 
  endofunctor
  $1+P: \slice I \to \slice I$, while $D$ is the
  set of such trees but with a marked leaf.  We construct these
  two objects as least fixpoints.  Put $Q = 1+P$; in  the internal language 
  we have
\[
Q(X_i \mid i \in I) = 
 \Big( \; \{ i \} + \sum_{a \in A_i} \prod_{b \in B_a} X_{sb} \mid
i \in I \; \Big) \,.
\]
Let $(C_i \mid i \in I)$ be the initial algebra for $Q$.
Its structure map is given by the family of isomorphisms
\begin{equation}
\label{equ:supC}
\mysup_{C_i} : \{ i \} + \sum_{a \in A_i} \prod_{b \in B_a} C_{sb} \isopil C_i
\,,
\end{equation}
meaning that a $Q$-tree is either a trivial tree (of some type $i\in I$)
or a one-node tree which is a term from $P$ (that is the choice of $a\in A_i$)
and whose incoming edges are labelled
by $Q$-trees (that is the map $k: B_a \to C_{sb}$).
We now define the polynomial endofunctor $R : \slice{C} \rightarrow \slice{C}$ 
by letting
\[
 R(X_c \mid c \in C) = \big( \widetilde{X}_c \mid c \in C \big) \,,
\]
where
\[
\widetilde{X}_c = 
\left\{
\begin{array}{ll}
\{ i \} &  \text{if } c = \sup(i) \, , \\[1ex]
 \sum_{b \in B_a} X_{kb} \quad & \text{if } c = \sup(a,k) \, .  
\end{array}
\right.
\]
This definition can be seen to be that of a polynomial functor using the
isomorphisms in~\eqref{equ:supC} and the extensivity of $\catE$. Let
$(D_c \mid c \in C)$ be the initial algebra for $R$.
Its structure maps consist of the following isomorphisms:
\[
\mysup_{D_{\mysup_C(i)}} : \{ i \}  \isopil D_{\mysup_C(i)} \, , \quad
\mysup_{D_{\mysup_C(a,h)}} : \sum_{b \in B_a} D_{hb} \isopil D_{\mysup_C(a,h)}
\,.
\]
The idea here is that a tree with a marked leaf is either a trivial tree,
with the unique leaf marked, or it is a pointed collection of trees, for which
the distinguished tree has a marked leaf.
We now define $u : D \rightarrow I$ recursively so that we have
\[
u(d) =
\left\{
\begin{array}{ll}
i & \text{if } d = \mysup_D(i) \, , \\
u(d') \qquad & \text{if } d = \mysup_D(b,d') \, .
\end{array}
\right.
\]
We have now constructed the polynomial in \eqref{equ:DC}, and we proceed to
verify that the associated polynomial functor is naturally isomorphic to $T$.
To prove this, it is sufficient to show that for every $X = (X_i \mid i \in I)$,
the object
\[
  \Big( \sum_{c \in C_i} \prod_{d \in D_c} X_{ud} \mid i \in I \Big)
\]
enjoys the same universal property that characterises $TX$, namely that
of being an initial algebra for the functor $X + P(-) : \slice{I} \rightarrow 
\slice{I}$. The required structure map is given by the following chain of 
isomorphisms:
\begin{eqnarray*} 
X_i + \sum_{a \in A_i} \prod_{b \in B_a} 
\sum_{c \in C_{sb}} \prod_{d \in D_c} X_{ud}
 & \cong & X_i + \sum_{a \in A_i} \sum_{k \in\!\! \underset{b \in 
 B_a}{\prod}\!\! C_{sb} } 
\prod_{b \in B_a} \prod_{d \in D_{kb}} X_{ud} \\
 & \cong & X_i + \sum_{(a,k) \in \!\!\underset{a \in A_i}{\sum} 
\underset{b \in B_a}{\prod}\!\! C_{sb}}
\prod_{(b,d) \in \!\!\underset{b \in B_a}{\sum} \!\! D_{kb}}  X_{ud}  \\
 & \cong & \sum_{c \in C_i} \prod_{d \in  D_c} X_{ud} 
\,.
\end{eqnarray*}
The initiality of the algebra follows by the initiality of $C$ and $D$
via lengthy, but not difficult, calculations. 

It remains to show that the unit and multiplication are cartesian. For
the unit $\eta : \Id{} \Rightarrow T$, we construct a diagram
\[
\xymatrix{
I \ar@{=}[d] & \ar@{=}[l] I\drpullback  \ar@{=}[r]  \ar[d]  & I \ar@{=}[r] 
\ar[d]^{e} &  I \ar@{=}[d]  \\
I & D \ar^{u}[l] \ar_{g}[r]     & C \ar[r]_{v} & I }
\]  
representing a cartesian strong natural transformation that coincides with 
$\eta$, modulo the 
isomorphism established above. 
For $i \in I$, we define $e_i : \{ i \} \rightarrow C_i$ by letting
$e_i(i) = \mysup_C(i)$. With this definition, we have an isomorphism 
$\{ i \} \cong D_{e_i(i)}$ for every $i \in I$, hence the middle
square is cartesian. We proceed analogously for the multiplication.
We will construct a diagram of the form
\[
\xymatrix{
I \ar@{=}[d] & \ar[l] F \drpullback \ar[r]  \ar[d]  & E \ar[r] \ar[d]^{m} &  I\ar@{=}[d]   \\
I  & D \ar^{u}[l] \ar_{g}[r]     & C \ar[r]_{v} & I}
\]  
where the top polynomial represents $T^2 : \slice{I} \rightarrow \slice{I}$
and the diagram represents the multiplication. Direct calculations with the
definition of substitution show that, for $i \in I$, we have
\[
E_i = \sum_{c \in C_i} \prod_{d \in D_c} C_{ud} \, , 
\]
and that, for $(c,k) \in E_i$, we have
\[
F_{(c,k)} = \sum_{d \in D_c} D_{kd} \, . 
\]
The family of maps $m_i : E_i \rightarrow C_i$ is defined recursively so that,
for $(c,k) \in E_i$, we have
\[
m_i(c,k) =
\left\{
\begin{array}{ll}
k(i) & \text{ if } c = \mysup_C(i) \\
\mysup(a, (\lambda b \in B_a) \; m_{sb}(hb, kb)) & \text{ if } c = \mysup_C(a,h)
\,.
\end{array}
\right.
\]
To check that the second clause is well-defined, observe that 
if $\mysup_C(a,h) \in C_i$ then, for $b \in B_a$, we have
$hb \in C_{sb}$. Furthermore we have 
\[
\prod_{d \in D_{\mysup(a,h)}} \! C_{ud} \; \cong \!
\prod_{ (b,d') \in \underset{b \in B_a}{\sum} D_{hb} } \!\! C_{u(b,d')} \; \cong
\prod_{b \in B_a} \prod_{d' \in D_{hb}} C_{u(d')}\,.
\]
Hence, for $b \in B_a$, we can regard $kb$ as an element of
\[
\prod_{d' \in D_{hb}} C_{u(d')}
\]
so that $(hb, kb) \in E_{sb}$, and therefore $m_{sb}(hb, kb) \in C_{sb}$, as
required. It is now easy to check that, for $(c,k) \in E_i$, we have
an isomorphism 
\[
  D_{m_i(c,k)} \cong F_{(c,k)} \, . 
\]
It remains to check that the natural transformation induced by the
diagram above is indeed the multiplication of the free monad on $P$.
This involves checking that its components satisfy the
condition that determines $\mu_X : T^2 X \rightarrow
TX$ uniquely. This is a lengthy calculation which we omit. 
\end{proof}

\begin{para} 
  To conclude this section, we derive from Theorem~\ref{thm:freemonad} a 
  stronger universal property of the free monad.
  Let $\PolyEnd_\catE$ denote the category whose objects are pairs $(I, P)$ consisting of an object $I \in \catE$
  and a polynomial endofunctor $P$ on $\slice{I}$, and whose morphisms from $(I,P)$
  to $(I', P')$ consist of a map $u : I' \rightarrow
  I$ in $\catE$ and a cartesian strong natural transformation
\begin{equation}
\label{equ:polyendmorphism}
\xymatrix{
\slice{I'} \ar[r]^{P'}  \ar[d]_{\ladj{u}} \ar@{}[dr]|{\ \Downarrow \, \phi} & \slice{I'} \ar[d]^{\ladj{u}} \\
\slice{I} \ar[r]_{P} & \slice{I} \,.}
\end{equation} 
The category $\PolyMnd_\catE$ of polynomial monads in $\catE$ is defined in a 
similar way: the objects are pairs $(I,T)$ consisting of an object $I \in \catE$
and a polynomial monad $T$ on $\slice{I}$. Maps from $(I,T)$ to $(I',T')$ are 
as in \eqref{equ:polyendmorphism}, but required now to satisfy the following 
monad map axioms:
\begin{equation}
\label{equ:genmonmor}
\xycenter{ \ladj{u} \ar[rr]^{\ladj{u} \comp \eta'} \ar@/_1pc/[drr]_{\eta \comp
\ladj{u}} & & \ladj{u} \comp T' \ar[d]^{\phi} \\
 & & T  \comp \ladj{u} } \qquad
 \xycenter{
 \ladj{u} \comp {T'}^2  \ar[r]^{\phi \comp T'}  \ar[d]_{\ladj{u} \mu'} & T \comp \ladj{u} 
 \comp T' \ar[r]^{T  \phi } & 
 {T}^2 \comp \ladj{u} \ar[d]^{\mu \ladj{u}} \\
 \ladj{u} \comp T' \ar[rr]_{\phi} & & T \comp \ladj{u} \,.}
 \end{equation}
Let us point out that the monad morphisms defined above are more special than
those that would arise by instantiating the notion of a monad morphism between
monads in a $2$-category, as defined in~\cite{StreetR:fortm}, to the $2$-category
$\PolyFunCart$: we allow only functors of the form $\ladj{u}
: \slice{I'} \rightarrow \slice{I}$, rather than arbitrary polynomial functors,
as vertical maps in the diagram~\eqref{equ:polyendmorphism}.  Note also that our
direction of $2$-cells are the oplax monad maps rather than the lax ones.
\end{para}

\begin{corollary}\label{thm:freemndglobal} 
  The forgetful functor $U : \PolyMnd_\catE \rightarrow \PolyEnd_\catE$ has a
  left adjoint.
\end{corollary}

\begin{proof}
  Both $\PolyEnd_\catE$ and $\PolyMnd_\catE$ are fibred over $\catE$ via the
  functors mapping an object $(I,-)$ to $I$, and $U$ is a fibred functor.
  Therefore, to define a left adjoint to $U$, it is sufficient to define left
  adjoints to the forgetful functors
  \[
  U_I : \PolyMnd_\catE(\slice{I}) \rightarrow \PolyEnd_\catE(\slice{I}) \, , 
  \]
  where $\PolyMnd_\catE(\slice{I})$ and $\PolyEnd_\catE(\slice{I})$ denote the
  fibre categories over
  $I \in \catE$.  But each $U_I$ has a left adjoint, sending $P$ to the free
  monad on $P$, cf.~Theorem~\ref{thm:freemonad}.
  It remains to observe that the canonical natural transformation 
  $\alpha : P \Rightarrow T$ (`insertion of generators') is 
  strong and cartesian. But we even have a polynomial representation of it:
  with notation as in the proof of Theorem~\ref{thm:freemonad},
  $\alpha$ is given by the diagram
  $$
  \xycenter{ 
  I\ar@{=}[d] & \ar[l] B  \drpullback\ar[r] \ar[d]& A \ar[d]^{\alpha_1} \ar[r]&I\ar@{=}[d] \\
  I & \ar[l] D  \ar[r] & C \ar[r]&I , \\
}
$$
cf.~\eqref{alpha} for the description of $\alpha$; the map $\alpha_1$ takes a
term in $A$ and interprets it as a tree with one node.  The map $B \to D$ is
described similarly but with a marked leaf.
\end{proof}

\begin{para}
  Observe that even if the forgetful functor $U : \PolyMnd_\catE \rightarrow
  \PolyEnd_\catE$ is fibred, its left adjoint is not.  The situation is
  analogous to the one represented in the diagram
  \[
  \xycenter{
  \Cat  \ar[rr] \ar[dr] & & \Grph  \ar[dl] \\
   & \Set & }
  \]
  where $\Cat$ is the category of small categories, and $\Grph$ is the category
  of directed, non-reflexive graphs.  The forgetful functor, mapping a category
  to its underlying graph, is a fibred functor, but its left adjoint, the free
  category functor, is not.
\end{para}

\section{$P$-spans, $P$-multicategories, and $P$-operads}
\label{sec:spans-ope}

\begin{para}
  Let $\Span_{\catE}$ denote the bicategory of spans in $\catE$, as introduced
  in \cite{BenabouJ:intb}.  Under the interpretation of spans as linear 
  polynomials (cf.~Example~\ref{ex:linear}),
  composition of spans (resp.~morphisms of spans) 
  agrees with composition of polynomials (resp.~morphisms of polynomials),
  so we can regard $\Span_{\catE}$ as a locally full sub-bicategory of
  $\PolyCart$, and view polynomials as a natural `non-linear'
  generalisation of spans.  
\end{para}
  
\begin{para}
  There is another notion of `non-linear' span, namely the $P$-spans of
  Burroni~\cite{BurroniA:tcat}, which is a notion relative to is a cartesian
  monad $P$.  This section is dedicated to a systematic comparison between the
  two notions, yielding (for a fixed polynomial monad $P$) an equivalence of
  framed bicategories between Burroni $P$-spans and polynomials over $P$ in the
  double-category sense.  We show how the comparison can be performed directly
  at the level of diagrams by means of some pullback constructions.  Considering
  monads in these categories, we find an equivalence between $P$-multicategories
  (also called coloured $P$-operads) and polynomial monads over $P$, in the
  double-category sense.

  In this section, strength plays no essential role: everything is cartesian
  relative to a fixed $P$, eventually assumed to be polynomial and hence strong,
  and for all the cartesian natural transformations into $P$ there is a unique
  way to equip the domain with a strength in such a way that the natural
  transformation becomes strong.
\end{para}

\begin{para}
  We first need to recall some material on $P$-spans and their extension.
  To avoid clutter, and to place ourselves in
  the natural level of generality, we work in a cartesian 
  closed category $\catC$, and consider a fixed
  cartesian endofunctor  $P : \catC \rightarrow \catC$.
  We shall later substitute $\catE/I$ for $\catC$,
  and assume that $P$ is a polynomial monad on $\catE/I$.
\end{para} 
 
\begin{para}
  By definition, a \myemph{$P$-span} is a diagram in $\catC$ of the form
  \begin{equation}
  \label{equ:pspan} 
  \xycenter{
  P(D) & \ar[l]_(.4){d}  N \ar[r]^{c} &  C \, , }
  \end{equation} 
A \myemph{morphism of $P$-spans} is a diagram like
  \begin{equation}
  \label{equ:pspanmorphism}
  \xycenter{
  P(D')  \ar[d]_{P(f)} & N' \ar[l]_(.4){d'} \ar[r]^c \ar[d]^g  & C' \ar[d]^h \\
  P(D)   & N \ar[r]_{c} \ar[l]^(.4){d} & C \, ,}
  \end{equation}
  We write $\PSpan$ for the category of $P$-spans and $P$-span morphisms in 
  $\catC$. 
\end{para} 

\begin{para}
  Let $\Cart_\catC$  denote the category whose objects are cartesian functors
between slices of $\catC$  and whose arrows are diagrams of the form
\begin{equation}
\label{psi}
\xycenter{
\catC/D' \ar[r]^{Q'}  \ar[d]_{\ladj{u}}  \ar@{}[dr]|{\ \Downarrow \, \psi}  
& 
\catC/C' \ar[d]^{\ladj{v}}  \\
\catC/D \ar[r]_{Q}    & \catC/C \,,}
\end{equation}
for $u : D' \rightarrow D$ and $v : C' \rightarrow C$ in $\catC$,
and $\psi$ a cartesian natural transformation. 
Under the identification $\catC=\catC/1$, we can consider $P$ as an object of 
$\Cart_\catC$, so it makes sense to consider the slice category
$\Cart_\catC/P$: its objects are the
cartesian functors $Q : \catC/D \rightarrow \catC/C$ equipped with a cartesian natural
transformation 
\begin{equation}
\label{equ:cartsq}
\xycenter{
\catC/D \ar[r]^{Q}  \ar[d]  \ar@{}[dr]|{\ \Downarrow \, \phi}
&  \catC/C \ar[d]    \\
\catC \ar[r]_{P}    & \catC \,.}
\end{equation} 

We now construct a functor $\ext : \PSpan_\catC \rightarrow \Cart_\catC/P$.  Its action on
objects is defined by mapping a $P$-span 
$PD \stackrel{d}\longleftarrow  N \stackrel{c}{\longrightarrow}  C$ to 
the diagram
\[
  \xymatrix{
  \catC/D  \ar[r]^(.44){P_{/D}}  \ar[d] & \catC/PD  \ar@{}[dr]|{\Downarrow }
  \ar[d] \ar[r]^{\pbk{d}} &
  \catC/N \ar[r]^{\ladj{c}} \ar[d] & \catC/C \ar[d] \\
  \catC \ar[r]_{P}  & \catC  \ar@{=}[r] & \catC \ar@{=}[r]  & \catC \,.}
  \]
Here $P_{/D} : \catC/D \rightarrow \catC/PD$  sends $f:X\to D$ to $Pf:  
PX \to PD$, and
the outer squares are commutative. The middle square is essentially
given by the counit of the adjunction $\ladj{d} \adjoint \pbk{d}$, and is
therefore a cartesian natural transformation.  More precisely, it is the
mate~\cite{KellyG:reve2c} of the commutative square
 \[
 \xymatrix{
 \catC/PD \ar[d]  & \catC/N \ar[l]_(.44){\ladj{d}}\ar[d] \\
 \catC  & \catC  \ar@{=}[l]  \,.
 }
  \]
The action of the functor $\ext : \PSpan_\catC \rightarrow
\Cart_\catC/P$ on morphisms is defined by mapping a diagram like
\eqref{equ:pspanmorphism} to the natural transformation
 \[
\xymatrixrowsep{35pt}
\xymatrixcolsep{35pt}
 \xymatrix{
 \catC/D' \ar[d]_{\ladj{f}} \ar[r]^{P_{/D'}} & \catC/PD'  \ar@{}[dr]|{\Downarrow} 
 \ar[d]_{\ladj{Pf}} 
 \ar[r]^{\pbk{d'}}  & \catC/N'  \ar[d]^{\ladj{g}} 
 \ar[r]^{\ladj{c'}}  & \catC/C' 
 \ar[d]^{\ladj{h}} \\
 \catC/D \ar[r]_{P_{/D}} & \catC/PD \ar[r]_{\pbk{d}} & \catC/N 
 \ar[r]_{\ladj{c}} & \catC/C \,,}
 \]
together with the structure maps to $P$.
The outer squares are just commutative and the middle square (again cartesian)
is the mate of the identity $2$-cell
\[
\xymatrixrowsep{35pt}
\xymatrixcolsep{35pt}
\xymatrix{
\catC/PD' \ar[d]_{\ladj{Pf}} & \catC/N'
\ar[l]_{\ladj{d'}} \ar[d]^{\ladj{g}} \\
\catC/PD  & \catC/N \ar[l]^{\ladj{d}}
\,.
}
\]
\end{para}

\begin{proposition}\label{thm:functofspan} 
  The functor $\ext : \PSpan_\catC \rightarrow \Cart_\catC/P$ is 
  an equivalence of categories.
\end{proposition} 
  
\begin{proof}
  The quasi-inverse is defined by mapping
  \[
  \xymatrix{
  \catC/D \ar[r]^{Q}  \ar[d]  \ar@{}[dr]|{\ \Downarrow \, \phi}
  & \catC/C \ar[d] \\
  \catC \ar[r]_{P} & \catC }
  \]
to the $P$-span 
\[
\xymatrix{ PD  & \ar[l]_{\phi_D} QD \ar[r]^{Q(1_D)} &  C \,.}
\]
The verification of the details is straightforward.
\end{proof}

\begin{para}
Given a cartesian natural transformation $\theta:P\Rightarrow P'$, there is a 
shape-change functor 
  \begin{eqnarray*}
    \PSpan_\catC & \longrightarrow & \PPSpan_\catC  \\
    {}[PD \leftarrow N \to C] & \longmapsto & [P'D\stackrel {\theta_D} \leftarrow PD \leftarrow 
    N \to C] \, .
  \end{eqnarray*}
  We also have the functor
    \begin{eqnarray*}
    \Cart_\catC/P & \longrightarrow & \Cart_\catC/P'  \\
    {}[Q\Rightarrow P] & \longmapsto & [Q\Rightarrow P \stackrel\theta 
    \Rightarrow P'] \, .
  \end{eqnarray*}
\end{para}

\begin{lemma}
  The equivalence $\ext$ of Proposition~\ref{thm:functofspan} is compatible with
  change of shape, in the sense that the following diagram commutes:
\[
\xymatrix{
\PSpan_\catC \ar[r]^{\ext} \ar[d] & \Cart_\catC/P \ar[d] \\
\PPSpan_\catC \ar[r]_{\ext} & \Cart_\catC/P' \, .}
\]
\end{lemma}
\begin{proof}
The claim amounts to checking
 \[
  \pbk{\theta_D}   \circ P'_{/D}  = P_{/D}
\]
which follows from the assumption that $\theta$ is cartesian.
\end{proof}

\begin{para}
We now assume that $P$ is a cartesian monad, so we have
two natural transformations $\eta : 1 \Rightarrow P$ and
$\mu : P\circ P \Rightarrow P$
at our disposal for shape-change.
As is well known \cite{LeinsterT:higohc}, this allows us to define horizontal
composition of $P$-spans: given composable $P$-spans 
\[
\xymatrix@dr{
N \ar[r]^c \ar[d]_d & C \\
PD } \qquad
\xymatrix@dr{
U \ar[r]^t \ar[d]_s & B \,, \\
PC }
\]
we define their composite $P$-span by applying $P$ to the first $P$-span, performing a pullback,
and using the multiplication map:
\begin{equation}\label{Pspan-comp}
  \xymatrixrowsep{36pt}
  \xymatrixcolsep{36pt}
  \vcenter{\hbox{
  \xymatrix @!0 {
      &     &            &   PN \times_{PC} U  \ar[dr] \ar[dl]  &      &  \\
       &    & PN  \ar[dl]_{Pd}  \ar[dr]^(.45){Pc}     &     &     U \ar[dl]_s 
	   \ar[dr]^t &  \\
&PPD \ar[ld]_{\mu_D}  &            & PC                            &      & B \\
PD  &&&&&
}}}
\end{equation}
Associativity of the composition law (up to coherent isomorphism)
depends on that fact that $P$ preserves pullbacks and that $\mu$ is cartesian.
It further follows from the fact that $\eta$ is cartesian that for each $D$
the $P$-span
\[
\xymatrix { 
PD & D \ar[l]_(.4){\eta_D} \ar[r]^{1_D} & D }
\]
is the identity $P$-spans for the composition law (up to coherent isomorphisms).
It is clear that these constructions are functorial in vertical maps between
$P$-spans, yielding altogether a double category of $P$-spans, denoted $\PSpan_\catC$:
the objects and vertical morphisms are those of $\catC$, the horizontal arrows
are the $P$-spans, and the squares are diagrams like (\ref{equ:pspanmorphism}).
\end{para}

\begin{para}
  We also have a double-category structure on $\Cart_\catC/P$: the horizontal
  composite of $Q\Rightarrow P$ with $R\Rightarrow P$ is $R\circ Q \Rightarrow P
  \circ P \Rightarrow P$, and the horizontal identity arrow is $\Id{}\Rightarrow
  P$.  Let us verify that the extension of a horizontal composite is isomorphic
  to the composite of the extensions: in the diagram
  \[
  \xymatrixrowsep{40pt}
  \xymatrixcolsep{48pt}
  \xymatrix @!0 {
  &&&\catC/C \ar[dr]^{P_{/C}}&&&\\
  &&\catC/N \ar[ur]^{\ladj{c}} \ar[dr]_{P_{/N}}&&\catC/PC \ar[dr]^{\pbk{s}} &&\\
  & \catC/PD \ar[ur]^{\pbk{d}} \ar[dr]^{P_{/PD}} &&\catC/PN \ar[ur]^{\ladj{(Pc)}} 
  \ar[dr] &\text{\tiny{B.C.}}& \catC/U \ar[dr]^{\ladj{t}} & \\
  \catC/D \ar[ru]^{P_{/D}} \ar[rr]_{(PP)_{/D}} && \catC/PPD \ar[ru]_{\pbk{(Pd)}} && 
  \catC/PN \times_{PC} U \ar[ru] && \catC/B
  }
  \]
  the top path is the composite of the extension functors, and the bottom path is
  the extension of the composite span.  The square marked B.C.~is the Beck-Chevalley
  isomorphism for the cartesian square~\eqref{Pspan-comp}, and the other squares, as well as the
  triangle, are clearly commutative.  The following proposition now follows from 
  Proposition~\ref{thm:functofspan}.
\end{para}

\begin{proposition}
  The functor $\ext: \PSpan_\catC \longrightarrow \Cart_\catC/P$ is an
  equivalence of double categories, in fact an equivalence of framed
  bicategories.
\end{proposition}
\noindent
We just owe to make explicit how the double category of $P$-spans is 
a framed bicategory: to each
vertical map $u: D' \to D$ we associate the $P$-span
$$
PD' \stackrel{\eta_{D'}}\longleftarrow D' \stackrel u \longrightarrow D \, .
$$
This is a left adjoint; its right adjoint is the $P$-span
$$
PD \stackrel{\eta_D \circ u} \longleftarrow D' \stackrel = \longrightarrow D'
$$
as follows by noting that their extensions are respectively
$\ladj{u}$ and $\pbk{u}$.  For this the important fact is that $\eta$ is 
cartesian.  With this observation it is clear that the  
equivalence is framed.

\begin{para}
  We now specialise to the case of interest, where $\catC=\catE/I$ and $P$ is a
  polynomial monad on $\catE/I$, represented by
  $$
  I \leftarrow B \to A \to I \,.
  $$
  Since now all the maps involved in the $P$-spans are over $I$,
  a $P$-span can be interpreted as a commutative diagram
  \[
  \xymatrix @dr {
    N \ar[d]_d \ar[r]^c & C \ar[d] \\
    PD \ar[r] & I \,.
  }
  \]
  If $C$ is an object of $\catC$, i.e.~a map in $\catE$ with codomain $I$, 
  we shall write $C$ also for its domain, and we have a natural identification of
  slices $\catC/C \simeq \catE/C$. That $P: \catE/I \to \catE/I$ is a polynomial monad, 
  means thanks to
  Lemma~\ref{thm:carttopoly}, that all objects in $\Cart_{\slice{I}}/P$ are polynomial again,
  so $\Cart_{\slice I}/P \iso \PolyCart/P$, the 
  category of polynomials cartesian over $P$ in the double-category sense.
  In conclusion:
\end{para}

\begin{proposition}\label{thm:PSpan=Poly/P}
  The functor $\ext:\PSpan_{\slice I} \to \PolyCart/P$ is an equivalence of
  framed bicategories.\qed
\end{proposition}

\begin{para}
  It is a natural
question whether there is a direct comparison between $P$-spans and polynomials over
$P$, without reference to their extensions.  This is indeed the case, as we now
proceed to establish, exploiting the framed structure.
  Given a polynomial over $P$, like
 \[
 \xymatrix{
Q: & D \ar[d]_u & M \ar[l] \ar[d] \ar[r] \drpullback & N \ar[d] \ar[r] & C \ar[d] \\
P: & I & B \ar[l]\ar[r] & A\ar[r] & I \,.}
\phantom{xxxx}
 \]
Consider the canonical factorisation of this morphism through the source lift of
$P$ along $u$ (cf.~\ref{sourcelift}):
\begin{equation} 
\label{QP-poly-span}
\xycenter{
Q: &D \ar@{=}[d] &  M  \ar[l] \ar[r]  \ar[d]   \drpullback & N   \ar[r]^c  
\ar[d]_d  & C  \ar[d] \\
\basechange{(u,\id{})} P :& D \ar@{=} [d] & \cdot  \drpullback  \ar[l] \ar[d] \ar[r]^f &  
PD  \ar[r] \ar[dd] & I  \ar@{=}[dd] \\
&D  \ar[d]_u  & \cdot \dlpullback  \ar[d] \ar[l] & & \\
P:&I  & \ar[l] B \ar[r] & A  \ar[r] & I \,. }
\phantom{xxxxxxxxx}
\end{equation}
Now we just read off the associated 
$P$-span:
\[
\xymatrix {
N \ar[r]^c \ar[d]_d & C \ar[d]\\
PD\ar[r] & I \,.}
\]
Conversely, given such a $P$-span, place it on top of the rightmost leg of
$P\circ \ladj{\alpha} = \basechange{(u,\id{})} P$ (the middle row of the
diagram, which depends only on $\alpha$ and $P$), and let $M$ be the pullback of
$N \to P(D)$ along the arrow labelled $f$.  It is easy to see that these
constructions are functorial,
yielding an equivalence of
hom categories $\PolyCart(D,C)/P \simeq \PSpan_{\slice I}(D,C)$.
\end{para}

\begin{example}
  Endo-$P$-spans $PC \leftarrow N \to C$, that is, polynomial endofunctors over
  $P$, are called \myemph{$C$-coloured $P$-collections}.  If furthermore $C=I$ we
  simply call them \myemph{$P$-collections}.  These are just polynomial
  endofunctors $Q : \catE/I \to \catE/I$ equipped with a cartesian natural
  transformation to $P$.  This category is itself a slice of $\catE$: it is
  easy to see that the functor
  \begin{eqnarray*}
    \PolyCart(I,I)/P & \longrightarrow & \catE/P1  \\
    Q & \longmapsto & [Q1\to P1]
  \end{eqnarray*}
  is an equivalence of categories.
\end{example}

\begin{para} 
  Burroni~\cite{BurroniA:tcat}, Leinster~\cite{LeinsterT:higohc}, and
  Hermida~\cite{HermidaC:repm} define $P$-multicategories (also called coloured
  $P$-operads) as monads in the bicategory of $P$-spans.  
  $P$-multicategories are also monads in the double category of $P$-spans --- this
  description also provides the $P$-multifunctors as (oplax) cartesian monad maps.
  $P$-multicategories based at the terminal object in $\catE/I$ are called
  $P$-operads.  If the base monad $P$ is a polynomial monad, the equivalence of
  Proposition~\ref{thm:PSpan=Poly/P} induces an equivalence of the categories of
  monads, as summarised in the corollary below.  
  
  In the classical example, $\catE$ is $\Set$ and $P$ is the free-monoid monad
  $M$ of Example~\ref{ex:freemonoid}.  In this case, $M$-multicategories are the
  classical multicategories of Lambek~\cite{LambekJ:dedsc2}, which are also
  called coloured nonsymmetric operads.  In the one-object case, $M$-operads are
  the plain (nonsymmetric) operads.  The other standard example is taking $P$ to
  be the identity monad on $\Set$.  Then $P$-multicategories are just small
  categories and $P$-operads are just monoids.  Hence small categories are
  essentially polynomial monads on some slice $\Set/C$ with an oplax cartesian
  double-categorical monad map to $\Id{}$, and monoids are essentially
  polynomial monads on $\Set$ with a cartesian monad map to $\Id{}$.  In
  summary, we have the following result.
\end{para}

\begin{corollary}
  There are natural equivalences of categories 
  \begin{xalignat*}{3}
    \PMultiCat & \simeq  \PolyMnd/P  &
    P\text{-}\defcat{Operad} & \simeq  \PolyMnd(1)/P \\
    \defcat{Multicat} & \simeq  \PolyMnd/M &
    \defcat{PlainOperad} & \simeq  \PolyMnd(1)/M \\
    \Cat & \simeq \PolyMnd/\Id{}  &
    \defcat{Monoid} & \simeq  \PolyMnd(1)/\Id{} \,.
  \end{xalignat*}
\end{corollary}

\begin{para}
  The double category of polynomials is very convenient for reasoning with
  $P$-multicategories.  The role of the base monad $P$ for $P$-multicategories
  is to specify a profile for the operations.  This involves specifying the
  shape of the input data, and it may also involve type constraints on input and
  output.  In the classical case of $P=M$, the fibres of $\Nat'\to\Nat$ 
  (Example~\ref{ex:freemonoid}) are
  finite ordinals, expressing the fact that inputs to an operation in a
  classical multicategory must be given as a finite list of objects.  In this
  case there are no type constraints imposed by $P$ on the operations.

For a more complicated example, let 
$P: \Set/\Nat\to \Set/\Nat$ be the free-plain-operad monad, which 
takes a collection (i.e.~an object in $\Set/\Nat$) 
and returns the free plain operad on it~\cite[p.135, p.145, p.155]{LeinsterT:higohc}.
This monad is polynomial (cf.~\cite{KockJ:polfo}): it is represented by
\[
\xymatrix {
\Nat & \operatorname{Tr}^{\bullet} \ar[l]_{s} \ar[r]^p  
& \operatorname{Tr} \ar[r]^t & \Nat \, ,
}
\]
where $\operatorname{Tr}$ denotes the
set of (isomorphism classes of) finite planar rooted trees, and $\operatorname{Tr}^{\bullet}$ denotes the
set of (isomorphism classes of) finite planar  rooted trees with a marked node.  The map $s$ returns the
number of input edges of the marked node; the map $p$ forgets the mark, and
$t$ returns the number of leaves. 
A $P$-multicategory $Q$ has a set of objects and a set of operations.
Each operation has its input slots organised as the
set of nodes of some planar rooted tree, since this is how the $p$-fibres look 
like.   Furthermore, there are type constraints: each
object of $Q$ must be typed in $\Nat$, via a number that we shall call the
\myemph{degree} of the object, and a compatibility
is required between the typing of operations and the typing of objects.
Namely, the degree of the output object of an operation must equal the total 
number of leaves of the tree whose nodes index the input, and the degree of the
object associated to a particular input slot must equal the number of incoming 
edges of the corresponding node in the tree.
All this is displayed with clarity by the fact
that $Q$ is given by a diagram
 \[
 \xymatrix{
Q: & D \ar[d]_\alpha & M \ar[l] \ar[d] \ar[r] \drpullback & N \ar[d]^\beta \ar[r] & 
 D\ar[d]^\alpha \\
P: &\Nat & \operatorname{Tr}^{\bullet} \ar[l] \ar[r] 
& \operatorname{Tr} \ar[r] & \Nat }
 \]
The typing of the operations is concisely given by the map $\beta$, and the
organisation of the inputs in terms of the fibres of the middle map of $P$ is
just the cartesian condition on the middle square.  The typing of objects is
encoded by $\alpha$ and the compatibility conditions, somewhat tedious to formulate
in prose, are nothing but commutativity of the outer squares.

Finite planar rooted trees can be seen as $M$-trees, where $M: \Set\to\Set$ is
the free-monoid monad (\ref{ex:freemonoid}).  Abstract trees, in turn, can be
seen as polynomial functors (\ref{ex:trees}): to a tree is
associated the polynomial functor
\[
A \leftarrow N' \to N \to A \, , 
\]
where $A$ is the set of edges, $N$ is the set of nodes, and $N'$ is the set
of nodes with a marked incoming edge.  Formally, an $M$-tree is a tree over $M$
in 
$\PolyCart$, that is to say a diagram
 \[
 \xymatrix{
A \ar[d] & N' \ar[l] \ar[d] \ar[r] \drpullback & N \ar[d] \ar[r] & 
 A\ar[d] \\
1 & \Nat' \ar[l] \ar[r] 
& \Nat \ar[r] & 1 \,.}
 \]
\end{para}

\providecommand{\bysame}{\leavevmode\hbox to3em{\hrulefill}\thinspace}
\providecommand{\MR}{\relax\ifhmode\unskip\space\fi MR }

\providecommand{\MRhref}[2]{%
  \href{http://www.ams.org/mathscinet-getitem?mr=#1}{#2}
}
\providecommand{\href}[2]{#2}


\end{document}